\documentclass[12pt]{article}
\usepackage{amssymb}
\usepackage{amsmath}
\usepackage{multirow}
\usepackage{hyperref}
\hypersetup{
  colorlinks   = true, 
  urlcolor     = blue, 
  linkcolor    = blue, 
  citecolor   = red 
}

\newtheorem{theorem}{Theorem}[section]
\newtheorem{proposition}[theorem]{Proposition}
\newtheorem{corollary}[theorem]{Corollary}
\newtheorem{lemma}[theorem]{Lemma}

\newcommand{\gby}[1]{\left\langle#1\right\rangle}
\newcommand{\abs}[1]{\left\vert#1\right\vert}
\DeclareMathOperator{\modulo}{mod}
\DeclareMathOperator{\tr}{Tr}
\DeclareMathOperator{\wt}{wt}
\newcommand{\set}[1]{\left\{#1\right\}}
\newcommand{\normordprod}[1]{\raisebox{-.7ex}{\ensuremath{\stackrel{\circ}{\circ}}}
											#1 \raisebox{-.7ex}{\ensuremath{\stackrel{\circ}{\circ}}}}

\DeclareMathOperator{\End}{End}

\allowdisplaybreaks[4]

\begin{document}

\title{Jacobi trace functions in the theory of vertex operator algebras}\author{Matthew Krauel\thanks{Supported by the Japan Society of the Promotion of Science (JSPS).} \\ {\small Institute of Mathematics, University of Tsukuba} \\
{\small and} \\ Geoffrey Mason\thanks{Supported by the National Science Foundation (NSF).} \\
{\small Department of Mathematics, UC Santa Cruz}}

\date{}
\maketitle

\abstract
\noindent
 We describe a type of $n$-point function associated to strongly regular vertex operator algebras 
$V$ and their irreducible modules.\ Transformation laws with respect to the Jacobi group are developed for $1$-point functions.\ For certain elements in $V$, the finite-dimensional space spanned by the $1$-point functions for the irreducible modules is shown to be a vector-valued weak Jacobi form.\ A decomposition of $1$-point functions for general elements is proved, and shows that such functions are typically quasi-Jacobi forms.\ Zhu-type recursion formulas are provided; they show how an $n$-point function can be written as a linear combination of $(n-1)$-point functions with coefficients that are quasi-Jacobi forms. \\

\maketitle

\section{Introduction}
\label{intro}
Let $V=(V, Y, \mathbf{1},\omega)$ be a vertex operator algebra (VOA) 
of central charge $\mathbf{c}$ with vacuum vector $\textbf{1}$ and Virasoro element $\omega$.\ 
For a state $v\in V$, the vertex operator determined by $v$ is generally denoted
\begin{eqnarray*}
Y(v, z) = \sum_{n\in \mathbb{Z}} v(n)z^{-n-1},
\end{eqnarray*}
where $v(n)$ is called the \emph{$n^{th}$ mode} of $v$. We also define operators $L(n)$ by 
\begin{eqnarray*}
Y(\omega, z) = \sum_{n\in \mathbb{Z}} L(n)z^{-n-2}
\end{eqnarray*}
for the vertex operator associated to $\omega$.\ The VOA $V$ carries the conformal grading into finite-dimensional subspaces $V = \oplus_{n\in \mathbb{Z}} V_n$,
where $V_n = \{v\in V\mid L(0)v=nv\}$.

In the present paper we deal exclusively with VOAs that are \emph{simple} and \emph{strongly regular}.\
Strong regularity of $V$ entails that it is \emph{rational, $C_2$-cofinite}, \emph{CFT-type}
(i.e., $V_0=\mathbb{C}\mathbf{1}$ and $V=\oplus_{n\geq 0} V_n$), and also $V_1$ consists of \emph{primary} states ($L(1)V_1=0$).\ These assumptions may be taken as the basic requirements for an axiomatic approach to rational conformal field theory.\ For a review of the theory of such vertex operator algebras, cf.\ \cite{Mason-Lattice}.\ A simple, strongly regular VOA satisfies the following additional properties (loc.\ cit.):
 \begin{enumerate}
 \item $V$ has only a finite number of (inequivalent) irreducible admissible modules, denoted by $M^1 ,\dots ,M^s$ (\cite{DLM-twisted,Zhu}).
 \item $V$ has a nonzero, invariant bilinear form $\gby{\cdot ,\cdot}\colon V\times V \to \mathbb{C}$.\ It is \emph{nondegenerate, symmetric}, and \emph{unique} when normalized so that $\gby{\textbf{1},\textbf{1}}=-1$ (\cite{FHL,Li-Bilinear}).
  \item $V_1$ is a \emph{reductive} Lie algebra with respect to the bracket $[u, v] = u(0)v$.\ Moreover,
  each homogeneous space of each irreducible module $M^r$  is a linearly reductive $V_1$-module.\ (This is proved in \cite{DM-Effective} for the adjoint module $V$.\ The more general case for $M^r$ can be proved similarly.)
 \item $V$ has a `square-bracket' grading such that
 \[
 V=\bigoplus_{n\geq 0} V_{[n]},
 \]
 where $V_{[n]}=\set{v\in V \mid L[0]v=nv}$ and $V_{[0]}=\mathbb{C}\textbf{1}$ (\cite{Zhu}).
 \end{enumerate}

 From now on, $\gby{\cdot ,\cdot}$ is the canonical invariant bilinear form
 normalized as in 3) above.

 Elements $h_1 ,\dots ,h_m \in V_1$ are said to satisfy \textit{Condition H} if
 \begin{enumerate}
 \item they are linearly independent,
 \item $h_1 (0), \dots ,h_m (0)$ are \emph{semisimple} operators on each module $M^1 ,\dots ,M^s$ with \emph{rational integer} eigenvalues,  \item $[h_i ,h_j]=0\ (1\leq i, j\leq m)$.
 \end{enumerate}
Thanks to the reductivity of $V_1$, any set of elements $\{h_i\}$ satisfying Condition H is contained in
 \emph{Cartan subalgebra} of $V_1$ (ie., a maximal abelian Lie subalgebra consisting of semisimple elements).\ Conversely, a Cartan subalgebra has a basis of elements that satisfy Condition H
 (\cite{Mason-Lattice}).\ It is easy to see that if the elements $\{h_i\}$ satisfy Condition H then
\begin{eqnarray*} 
h_i (n) h_j =\delta_{n,1} \gby{h_i ,h_j}\textbf{1}\ (1\leq i, j\leq m).
\end{eqnarray*}

Until further notice, fix $h_1, \hdots, h_m$ satisfying Condition H.\ Introduce the symmetric matrix $G=(\gby{h_i ,h_j})$ and let $G[\underline{\alpha}]$ denote $\underline{\alpha}^tG \underline{\alpha}$ ($t$ denotes transpose) for an $m$-rowed column vector $\underline{\alpha}$.\
For elements $v_1 ,\dots ,v_n \in V$ we  consider \emph{$n$-point functions} of the form
\begin{equation}
 \begin{aligned}
  &F_{r}(v_1, \hdots, v_n; z_1, \hdots, z_m, \tau)\\
  &\hspace{10mm}=\tr_{M^r} Y(q_{1}^{L(0)} v_1 ,q_{1}) \cdots Y(q_{n}^{L(0)} v_n ,q_n) \zeta_1^{h_1 (0)}\cdots \zeta_m^{h_m (0)} q^{L(0)-\mathbf{c}/24},
  \end{aligned}
  \label{intro1}
 \end{equation}
where $\mathbf{c}$ is the central charge of $V$ and we always take $q_k =e^{w_k}, \zeta_l =e^{2\pi i z_l}, q=e^{2\pi i\tau}$
with $w_k, z_l \in \mathbb{C}\ (1\leq r\leq s, 1\leq k \leq n, 1 \leq l \leq m)$ and $\tau\in \mathbb{H}$
(the complex upper half-plane).\ It is convenient to abbreviate tuples such as
$(h_1, \hdots, h_m)$ by $\underline{h}$.\ Thus the important special case of a $1$-point function with homogeneous element $v\in V_k$, for example, reduces to
 \begin{equation}
 J_{r,\underline{h}}(v;\tau ,\underline{z}):=\tr_{M^r} o(v) \zeta_1^{h_1 (0)}\cdots \zeta_m^{h_m (0)} q^{L(0)-\mathbf{c}/24}, \label{intro2}
 \end{equation}
which (formally) can be written
 \begin{equation}
 J_{r,\underline{h}}(v;\tau ,\underline{z})=q^{\lambda_r -\mathbf{c}/24} \sum_{\ell \geq 0} \sum_{t_1 ,\dots ,t_m \in \mathbb{Z}} c(\ell ,t_1 ,\dots ,t_m) \zeta_1^{t_1} \cdots \zeta_m^{t_m} q^\ell ,
 \label{FourierExpansion}
 \end{equation}
 where $\lambda_r$ is the conformal weight of $M^r$.\ (In  case the trace is over a space $W$ which is not $M^r$, we will denote (\ref{intro2}) by $J_{W,\underline{h}}(v;\tau ,\underline{z})$.)
 
  For example, if $v=\mathbf{1}$ and $h_1, \hdots, h_m$ is a \emph{basis} of a Cartan subalgebra $H$ of $V_1$, then $J_{r, \underline{h}}(\mathbf{1}; \tau, \underline{z})$ (the $0$-point function) determines the multiplicities of the $H$-weights of $M^r$ considered as $V_1$-module, and thereby the decomposition of $M^r$ into irreducible $V_1$-modules.
 
The main purpose of the present paper is to establish transformation laws for one-point functions
 with respect to the \emph{Jacobi group} $SL_2(\mathbb{Z}) \ltimes (\mathbb{Z}\oplus\mathbb{Z})^m$.\ It transpires that this naturally breaks down into two cases, depending on whether the equality
 $h_j (n)v=0\ (1\leq j\leq m, n\geq 0$) holds or not.\ In the first case we have the following theorem.
 
  \begin{theorem} \label{TheoremAlpha}
 Let $V$ be a simple, strongly regular VOA.\  For any $v\in V$, there are finitely many integers $t$ such that the function $J_{r,\underline{h}}(v;\tau ,\underline{z})$ converges on every closed subset of $\{(\tau ,z_1 ,\dots ,z_m)\in \mathbb{H} \times \mathbb{C}^m \mid z_i \not \in \frac{1}{t} \left(\mathbb{Z} +\mathbb{Z}\tau\right), 1\leq i\leq m\}$ with Fourier expansion (\ref{FourierExpansion}).\ If $v\in V_{[k]}$ satisfies $h_j (n)v=0$ for all $1\leq j \leq m$ and $n\geq 0$, then $J_{r,\underline{h}}(v;\tau ,\underline{z})$ satisfies the following functional equations:
 \begin{enumerate}
 \item For all $\gamma =\left( \begin{smallmatrix} a&b\\c&d \end{smallmatrix}\right) \in \text{SL}_2 (\mathbb{Z})$,
\begin{align}\label{MainTheorem1}
    & J_{r,\underline{h}}\left(v; \frac{a\tau +b}{c\tau +d}
 ,\frac{\underline{z}}{c\tau +d}\right)=(c\tau +d)^k \exp \left(\pi i \frac{cG[\underline{z}]}{c\tau +d}\right) \sum_{\ell =1}^s A_{r,\gamma}^\ell J_{\ell ,\underline{h}}(v;\tau ,\underline{z}),
 \end{align}
with scalars $A_{r,\gamma}^\ell$ depending only on $\gamma$.
 \item For all $[\underline{\lambda},\underline{\mu}]\in \mathbb{Z}^m \times \mathbb{Z}^m$ there is a permutation $r\mapsto r'$, $r' \in \set{1, \dots ,s}$, such that
 \begin{equation}
 J_{r,\underline{h}}\left(v;\tau ,\underline{z}+\underline{\lambda}\tau +\underline{\mu} \right)=\exp \left(-\pi i (G[\underline{\lambda}]\tau +2\underline{z}^t G\underline{\lambda})\right) J_{r' ,\underline{h}}(v;\tau ,\underline{z}).
 \label{MainTheorem2}
 \end{equation}
 \end{enumerate}
 \end{theorem}  

Essentially, this says that the vector of $1$-point functions $(J_{1,\underline{h}}, \dots ,J_{s,\underline{h}})^t$
  is a \emph{vector-valued weak Jacobi form of weight $k$ and index $G/2$}.\ For example, if $V$ is 
  \emph{holomorphic} (i.e., it has a unique irreducible module), then $J_{V,\underline{h}}$ is a weak Jacobi form of weight $k$ and index $G/2$ (generally with a character $\chi$ of $SL_2(\mathbb{Z})$,
  which is trivial if $24|\mathbf{c}$).
  
  The statement of convergence can be refined in a number of cases. In particular, Heluani and Van Ekeren have recently \cite{HE-characters} introduced the idea to use another set of quasi-Jacobi forms to address this issue in the case of $N=1$ SUSY vertex algebras. In this setting, they first prove their trace functions are conformal blocks, and are then able to utilize this other set of quasi-Jacobi forms to prove convergence in the $m=1$ case on the stronger domain consisting of closed subsets of $\{(\tau ,z) \mid z\not \in \mathbb{Z} +\mathbb{Z}\tau\}$. It appears possible and of interest to extend these ideas to establish a similar domain of convergence in the case of general strongly regular VOAs.
 
When $v\in V_{[k]}$ fails to satisfy $h_j (n)v =0$ for some $1\leq j\leq m$ or $n\geq 0$, 
the $1$-point functions (\ref{intro2}) do not necessarily satisfy (\ref{MainTheorem1}) and (\ref{MainTheorem2}).\ To describe the transformation laws in this case, let us fix for now a Cartan subalgebra $H\subseteq V_1$, say of dimension $d$, together with an \emph{orthogonal basis}
$\left\{u_j\right\}$ of $H$.\  It suffices to take $v\in V_{[k]}$ in the form
\begin{equation}\label{introeq2}
    v=u_1 [-m_{1,1}]^{\ell_{1,1}} \cdots u_1 [-m_{1,\nu_1}]^{\ell_{1,\nu_1}}
    \cdots  u_d [-m_{d,1}]^{\ell_{d,1}}\cdots u_d
    [-m_{d,\nu_d}]^{\ell_{d,\nu_d}} w \hspace{-1ex}
 \end{equation}
 for nonnegative integers $\ell_{x,y} ,m_{x,y}$ ($1\leq x \leq d, 1\leq y\leq \nu_d$), and $w$ in the commutant $\Omega(0):= C_V(M_H)$ of the Heisenberg subVOA $M_H\subseteq V$ generated by $H$.\ There is a decomposition
  (\cite{DM-Effective,DM-Integrability,Mason-Lattice})
 \begin{align}
 M^r&=\bigoplus_{t=1}^\delta \bigoplus_{\beta \in \Lambda} M_H (\beta +\gamma_t) \otimes \Omega_{r} (\gamma_t) \label{Decomp21}
 \end{align}
of $M^r$ into irreducible $M_H\otimes \Omega(0)$-modules.\ In particular, each $\Omega_r(\gamma_t)$ is a certain irreducible $\Omega(0)$-module.\ Here, $\Lambda\subseteq P\subseteq H$ where $\Lambda, P$ are additive subgroups of $H$ of rank $d$, $\Lambda$ is a positive-definite even lattice with respect to $\langle \ , \ \rangle$, and $\set{\gamma_t }$ are coset representatives of $P/\Lambda$.\ (See \cite{Mason-Lattice} and Section 5 below for further details.)\ We then have the following theorem.

 \begin{theorem}\label{TheoremBeta}
 Let $V$ be a simple, strongly regular vertex operator algebra with orthogonal basis $\left\{u_j\right\}$ of $H$ and $v\in V_{[k]}$ as in (\ref{introeq2}).\
 Then
\begin{align}
J_{r,\underline{h}}(v,\tau ,\underline{z})&= \frac{1}{\eta (\tau)^{d}} \sum_{t=1}^\delta J_{\Omega_{r} (\gamma_t),\underline{h}}(w;\tau ,\underline{z}) \sum_{\underline{i}} f^{\underline{i}} (\tau)
\Psi_{t, \underline{h}}(k_{\underline{i}}, \tau, \underline{z}), \notag \label{TheoremBetaEq1}
\end{align}
where $\underline{i}=(i_{1,1}, \dots,i_{1,\nu_d},\dots ,i_{d,1},\dots ,
i_{d,\nu_d})\!\in\! \mathbb{Z}^{d \nu_d}$ with $0\!\leq\! i_{x,y} \!\leq\! \lfloor
\ell_{x,y}\rfloor$, $f^{\underline{i}} (\tau)$ is a quasi-modular form of weight 
$2\sum_{x,y} i_{x,y}, k_{\underline{i}}=\sum_{x,y} (\ell_{x,y}-2i_{x,y})$, and $\Psi_{t ,\underline{h}}(k_{\underline{i}} ,\tau ,\underline{z})$ is a linear combinations of functions of the form
\begin{equation}
\sum_{\alpha \in \Lambda +\gamma_t} \gby{a ,\alpha}^{k_{\underline{i}}} q^{\gby{\alpha ,\alpha}/2} \zeta_1^{\gby{\alpha ,h_1}} \cdots \zeta_m^{\gby{\alpha ,h_m}} \label{ThetaIntro}
\end{equation}
for various $a \in H$.\ (For precise definitions of these functions, see
Section~5.)
\end{theorem} 

The functions (\ref{ThetaIntro}) and their transformation laws with respect to the Jacobi group
 are discussed in \cite{KrauelI}.\ In the case $\gby{a, h_j}=0$ for all $1\leq j\leq m$, they are Jacobi forms on $\Gamma_0 (N)$ of weight $f+k_{\underline{i}}$ and index $G/2$, where $N$ is the level of $\Lambda$, $\Gamma_0 (N)\subseteq \text{SL}_2 (\mathbb{Z})$ is defined by
 \[
 \Gamma_0 (N)=\left \{\left( \begin{smallmatrix} a&b \\c&d \end{smallmatrix}\right) \in \text{SL}_2 (\mathbb{Z}) \mid c \equiv 0\ (\modulo N) \right \},
 \]
and we suppose the quadratic form has rank $2f$.\ Otherwise, (\ref{ThetaIntro}) are quasi-Jacobi forms on $\Gamma_0 (N)$ of the same weight and index. In either case, 
 \[
 \sum_{\underline{i}} f^{\underline{i}} (\tau) \Psi_{t, \underline{h}}(k_{\underline{i}}, \tau, \underline{z})
 \]
 is a quasi-Jacobi form on $\Gamma_0 (N)$ of weight $f+ \sum_{j} \ell_j$ and index $G/2$.

 Theorem \ref{TheoremBeta} reduces the computation of $J_{r, \underline{h}}(v, \tau, \underline{z})$ to a similar computation involving only the commutant $\Omega(0)$ and its irreducible modules.\ It is a standard conjecture that, under the assumption that $V$ is strongly regular, 
$\Omega (0)$ is also strongly regular.\ Assuming this to be true (it is known in many cases), Theorems \ref{TheoremAlpha} and \ref{TheoremBeta} provide explicit transformation laws for the functions $J_{r,\underline{h}}(v;\tau ,\underline{z})$ for any homogeneous $v\in V$. 

 The literature dealing with (weak) Jacobi forms in the context of affine algebras and related areas is quite extensive, whereas the
 theory for general vertex operator algebras that we develop here has few precedents.\ In \cite{DLMa-elliptic} some of the theory is developed for lattice VOAs, and  \cite{Kac-infinite} deals with the case of  highest weight integrable representations for affine Kac-Moody Lie algebras.\ Weak Jacobi forms arise as elliptic genera in various contexts,
 e.g., from models of $N=2$ super conformal field theories discussed in \cite{KYY}.\ In particular, a generic approach is developed to calculate the relevant transformation properties for $N=2$ Neveu-Schwarz models and the elliptic genus for the $N=2$ Landau-Ginzburg models are calculated.\ Libgober also discusses elliptic genera in \cite{Libgober-Elliptic}, showing in the Calabi Yau case that the elliptic genus is a weak Jacobi form, while in other cases it 
 lies in the space of quasi-Jacobi forms.

 Weak Jacobi forms and quasi-Jacobi forms also appear, at least implicitly, in the study of $n$-point recursion formulas in \cite{DLM-orbifold} and \cite{MTZ}.\ Gaberdiel and Keller \cite{GK-differential} discuss these functions further in the $N=2$ Neveu-Schwarz model, developing some transformation properties while also establishing differential operators which arise in superVOAs that preserve the weak Jacobi form property of the elliptic genus.\ Recent work of Heluani and Van Ekeren \cite{HE-characters} considers certain supercurves and the vertex (operator) algebras ($N_W =1$ SUSY vertex algebras) that produce vector bundles over these supercurves.\ In this setting, they show that certain functions analogous to the ones studied here give rise to superconformal blocks on a moduli space of elliptic supercurves.\ Their work establishing convergence of the functions they consider inspired us to revisit the convergence of the functions that we deal with here, where a previous draft of this paper contained an incomplete proof.\ As Heluani and Van Ekeren explain, their work can be regarded as an algebro-geometric approach to such problems.

 Finally, we note that the special case of the partition function with \emph{only one} elliptic variable, i.e., $v=\textbf{1}$ and $m=1$ in previous notation, appears in \cite{KrauelMasonI}.\ It is our hope that the results of the present paper may, in particular, foster closer ties between vertex operator algebras and elliptic genera.
 
 The paper is organized as follows.\ In Section \ref{Section-AutomorphicForms} we discuss the various kinds of modular-type functions that we need, including (matrix) Jacobi and
 quasi-Jacobi forms, and `twisted' Weierstrass and Eisenstein series.\ The latter functions appear as coefficients in the recursion formula, expressing $n$-point functions  (\ref{intro1}) as a sum of $(n-1)$-point functions, which is proved in the short Section \ref{Section-RecursionFormula}, following 
\cite{MTZ}.\ This result reduces the study of $n$-point functions to the case of
$1$-point functions.\ Our main results, Theorems  \ref{TheoremAlpha} and \ref{TheoremBeta}, are proved in Sections \ref{Section-TheoremAlpha} and \ref{Section-TheoremBeta} respectively.

\section{Automorphic forms\label{Section-AutomorphicForms}}

\subsection{Jacobi and quasi-Jacobi forms\label{Subsection-JacobiForms}} 
 
 Let $\text{Mer}_{\mathbb{H}\times \mathbb{C}^n}$ denote the space of meromorphic functions 
 on $\mathbb{H}\times \mathbb{C}^n$, and $F$ be a real symmetric positive-definite $n\times n$ matrix. We say a meromorphic function $\phi$ on $\mathbb{H}\times \mathbb{C}^n$ is a \textit{meromorphic Jacobi form of weight $k$, index $F$, and character $\chi$ $(\chi \colon \Gamma_1 \to \mathbb{C}^*)$} on a subgroup $\Gamma_1$ of $\text{SL}_2(\mathbb{Z})$ if for some $\ell_0 \in \mathbb{Q}$, $\phi$ has an expansion of the form
 \begin{equation}
 \phi (\tau ,\underline{z})= \sum_{\begin{subarray}{c}\underline{r}\in \mathbb{Z}^n, \ell \in
     \mathbb{Q},\\ 4\ell -F^{-1}[\underline{r}] \geq 0\end{subarray}} c(\ell ,\underline{r}) q^\ell \exp \left( 2\pi i (\underline{z}^t \underline{r}) \right), \label{JacobiExpansion1}
 \end{equation}
 where $q=e^{2\pi i \tau} \ (\tau \in \mathbb{H})$, $\ell \geq \ell_0$, $c(\ell ,\underline{r})$ are scalars, and for all $\gamma =\left( \begin{smallmatrix} a& b\\c &d \end{smallmatrix} \right) \in \Gamma_1$ and $(\underline{\lambda},\underline{\mu}) \in \mathbb{Z}^n \times \mathbb{Z}^n$ we have
 \begin{equation}
 \phi \left(\frac{a\tau +b}{c\tau +d}, \frac{\underline{z}}{c\tau +d} \right ) = \chi (\gamma) (c\tau +d)^{k} \exp \left( 2\pi i \frac{cF[\underline{z}]}{c\tau +d} \right) \phi (\tau ,\underline{z}), \notag 
 \end{equation}
 and
 \begin{equation}
  \phi (\tau ,\underline{z}+\underline{\lambda}\tau +\underline{\mu}) = \exp \left(-2\pi i(\tau F[\underline{\lambda}]+2\underline{z}^t F \underline{\lambda}) \right) \phi (\tau ,\underline{z}). \notag 
 \end{equation}
 In the case $\ell_0 \geq 0$, $\phi$ is \textit{holomorphic}. Throughout this paper we take the term Jacobi form to mean holomorphic Jacobi form. When the condition $4\ell -F^{-1}[\underline{r}] \geq 0$ in (\ref{JacobiExpansion1}) is replaced with $\ell \geq 0$, we call $\phi (\tau ,\underline{z})$ a \textit{weak Jacobi form of weight $k$ and index $F$}. (See \cite{EZ} for a detailed study of such functions when $n=1$ and \cite{Skor-critical} for a discussion of the general case.)
 
 The function $\phi$ is a \textit{quasi-Jacobi form} of weight $k$ and index $F$ on $\Gamma_1$ if for each $\tau \in \mathbb{H}$, $\underline{z}\in \mathbb{C}^n$, $\gamma =\left(\begin{smallmatrix}a&b \\c&d \end{smallmatrix}\right) \in \Gamma_1$, and $[\underline{\lambda},\underline{\mu}]\in \mathbb{Z}^n \times \mathbb{Z}^n$, we have 
 \begin{enumerate}
 \item $(c\tau +d)^{-k} e^{-2\pi i \frac{cF[\underline{z}]}{c\tau +d}} \phi \left( \frac{a\tau +b}{c\tau +d}, \frac{\underline{z}}{c\tau +d} \right) \in \text{Mer}_{\mathbb{H}\times \mathbb{C}^n} \left[\frac{cz_1}{c\tau +d}, \dots ,\frac{cz_n}{c\tau +d}, \frac{c}{c\tau +d}\right]$ with coefficients dependent only on $\phi$, and
 \item $e^{2\pi i(\tau F[\underline{\lambda}]+2\underline{z}^t F
     \underline{\lambda})} \phi (\tau ,\underline{z}+\underline{\lambda}\tau
     +\underline{\mu})\! \in\! \text{Mer}_{\mathbb{H}\times \mathbb{C}^n}[\lambda_1 ,\dots ,\lambda_n]$ with coefficients dependent only on $\phi$.
 \end{enumerate}
  In other words, there are meromorphic functions $S_{i_1 ,\dots ,i_n ,j}(\phi)$ and $T_{i_1 ,\dots ,i_n} (\phi)$ on $\mathbb{H} \times \mathbb{C}^n$ determined only by $\phi$, and $s_1 ,\dots ,s_n ,t \in \mathbb{N}$ such that
\begin{equation}
  \begin{aligned}
  &(c\tau +d)^{-k} \exp \left(-2\pi i \frac{cF[\underline{z}]}{c\tau +d}\right) \phi \left( \frac{a\tau +b}{c\tau +d}, \frac{\underline{z}}{c\tau +d} \right) \\
       &\hspace{5mm}=\sum_{\begin{subarray}{c}i_1 \leq s_1,\dots
  ,i_n \leq s_n\\ j\leq t\end{subarray}} S_{i_1 ,\dots ,i_n,j} (\phi)(\tau ,\underline{z})\left(\frac{cz_1}{c\tau +d}\right)^{i_1}\cdots \left(\frac{cz_n}{c\tau +d}\right)^{i_n} \left( \frac{c}{c\tau +d}\right)^j
 \end{aligned}\label{QuasiDefn1}
 \end{equation}
 and
 \begin{equation}
 \begin{aligned}
 &e^{2\pi i(\tau F[\underline{\lambda}]+2\underline{z}^t F \underline{\lambda})} \phi (\tau ,\underline{z}+\underline{\lambda}\tau +\underline{\mu}) \\
&\hspace{5mm}=\sum_{i_1 \leq s_1 ,\dots ,i_n \leq s_n} T_{i_1 ,\dots ,i_n} (\phi)(\tau ,\underline{z})\lambda_1^{i_1}\cdots \lambda_n^{i_n} . 
\end{aligned}\label{QuasiDefn2}
 \end{equation}
 If $\phi \not =0$, we take $S_{s_1 ,\dots ,s_n ,t}(\phi) \not =0$ and $T_{s_1 ,\dots ,s_n} (\phi) \not =0$, and say $\phi$ is a quasi-Jacobi form of \textit{depth} $(s_1 ,\dots ,s_n ,t)$. In the case $n=1$ and $F=0$, this definition of a quasi-Jacobi form reduces to that in \cite{Libgober-Elliptic}. (See also Definition $3.10$ in \cite{Kaw} for another definition of quasi-Jacobi form.)

 Let $\mathcal{Q}_n$ denote the space of quasi-Jacobi forms on $\mathbb{H}\times \mathbb{C}^n$. Straightforward calculations establish the following well-known lemma.
  
\begin{lemma} \label{QuasiClosedPartials}
 The space of quasi-Jacobi forms $\mathcal{Q}_n$ is closed under multiplication by $E_2 (\tau)$ and partial derivatives $\frac{d}{d\tau}$ and $\frac{d}{dz_j}$, $1\leq j\leq n$.\ (Nb.\ such operations change the weight and depth.)
\end{lemma} 
 
 \subsection{Twisted elliptic functions\label{Subsection-TwistedEllipticFunctions}}
 
  For $w\in \mathbb{C}$, $\underline{z} \in \mathbb{C}^n$, and $\tau \in \mathbb{H}$ such that $\abs{q}<\abs{e^{w}}<1$ and $\zeta_{z_1 +\cdots +z_n}\not =1$, we define the \textit{`twisted' Weierstrass functions} $\tilde{P}_k (w,\underline{z},\tau)$ by
  \begin{equation}
  \tilde{P}_k (w,\underline{z},\tau):= \frac{1}{(k-1)!} \sum_{\ell \in \mathbb{Z}} \frac{\ell^{k-1} q_w^\ell}{1-\zeta_1^{-1} \cdots \zeta_n^{-1} q^\ell},
  \label{twistedPk}
  \end{equation}
  where $q=e^{2\pi i \tau}$, $q_w =e^{w}$, $\zeta_j =e^{2\pi i z_j}$. When $z=z_1 +\cdots +z_n$ and we set $\zeta =e^{2\pi i z}$, the functions $\tilde{P}_k (w,z,\tau)=\tilde{P}_k (w,\underline{z},\tau)$ are the same as the functions $(-1)^k P_k \left[\begin{smallmatrix} \zeta \\ 1 \end{smallmatrix}\right] (w,\tau)$ in \cite{MTZ} (where one can find more details), $P_k (1,\zeta^{-1} ,2\pi iw,\tau)$ in \cite{DLM-orbifold}, and $(2\pi i)^{-k} \hat{\mathcal{P}}_{k}(q_w ,q,\zeta)$ in \cite{GK-differential}. We will also consider functions of the form (\ref{twistedPk}) when $\zeta_1 \cdots \zeta_n =1$. In this case, the sum is to exclude the term $\ell =0$ and the functions are simply the classical (or `untwisted') Weierstrass functions.

 Writing (\ref{twistedPk}) as
 \begin{align}
 \tilde{P}_k (w,\underline{z},\tau)&=\frac{1}{(k-1)!}\sum_{\ell=1}^\infty \left( \frac{\ell^{k-1}q_w^\ell}{1- \zeta_1^{-1} \cdots \zeta_n^{-1} q^\ell}+\frac{(-1)^k \ell^{k-1} q_w^{-\ell}q^\ell \zeta_1 \cdots \zeta_n}{1- \zeta_1 \cdots \zeta_n q^\ell}\right) \notag \\
 &\quad +\delta_{k,1} \frac{1}{1-\zeta_1^{-1} \cdots \zeta_n^{-1}}, \notag 
 \end{align}
 it can be shown that the functions $\tilde{P}_k (w,\underline{z},\tau)$ converge for $\vert q \vert < \vert q_w \vert <1$ (see also \cite{GK-differential,MTZ}).

 Define the functions $\tilde{G}_k (\tau ,\underline{z})$ by
  \begin{align}
  &\tilde{G}_{2k} (\tau ,\underline{z}):= 2\xi (2k) +\frac{(2\pi i )^{2k}}{(2k-1)!} \sum_{\ell =1}^\infty \left( \frac{\ell^{2k-1}q^\ell \zeta_1^{-1}\cdots \zeta_n^{-1}}{1-q^\ell \zeta_1^{-1}\cdots \zeta_n^{-1}} +\frac{\ell^{2k-1}q^\ell \zeta_1 \cdots \zeta_n}{1- q^\ell \zeta_1 \cdots \zeta_n} \right), \notag \\
  &\tilde{G}_{2k+1} (\tau ,\underline{z}):= \frac{(2\pi i )^{2k+1}}{(2k)!} \sum_{\ell =1}^\infty \left( \frac{\ell^{2k}q^\ell \zeta_1^{-1} \cdots \zeta_n^{-1}}{1-q^\ell \zeta_1^{-1} \cdots \zeta_n^{-1}} -\frac{ \ell^{2k}q^\ell \zeta_1 \cdots \zeta_n}{1- q^\ell \zeta_1 \cdots \zeta_n} \right), \notag \\
  &\tilde{G}_1 (\tau ,\underline{z}):= (2\pi i )\sum_{\ell=1}^\infty
      \left(\frac{q^\ell \zeta_1^{-1} \cdots \zeta_n^{-1}}{1- q^\ell
      \zeta_1^{-1} \cdots \zeta_n^{-1}}+\frac{q^\ell \zeta_1 \cdots \zeta_n}{1-
  q^\ell \zeta_1 \cdots \zeta_n}\right)+\frac{2\pi i}{1-\zeta_1^{-1} \cdots \zeta_n^{-1}} -\pi i, \notag
  \end{align}
  where $\xi (2k)=\sum_{n=1}^\infty \frac{1}{n^{2k}}$, and set
  \begin{equation}
  \tilde{E}_{m}(\tau ,\underline{z}):= \frac{1}{(2\pi i)^m} \tilde{G}_{m} (\tau ,\underline{z}). \notag
  \end{equation}
  The functions $\tilde E_m$ with one complex variable $z$ have been called `twisted Eisenstein series' in \cite{DLM-orbifold,GK-differential,MTZ}. The additional complex variables considered here do not add much difficulty as most calculations reduce to the single complex variable case by noting
  \[
  \tilde{E}_m (\tau ,\underline{z})=\tilde{E}_{m} (\tau ,z_1 +\cdots +z_n).
  \]
 
\begin{lemma}\label{quasi-JacobiEisenlemma}
 For $m\geq 1$, the functions $\tilde{E}_m (\tau ,\underline{z})$ are quasi-Jacobi forms of weight $m$ and index $0$.
\end{lemma}

\paragraph{Proof}
We first take the $n=1$ ($\underline{z}=z$) case and show that $\tilde E_m$ satisfies (\ref{QuasiDefn1}) for the matrices $S=\left(\begin{smallmatrix}0&-1\\1&0 \end{smallmatrix} \right)$ and $T=\left(\begin{smallmatrix}1&1\\0&1 \end{smallmatrix} \right)$. The result follows from a transformation discussed in \cite{GK-differential}. In particular, it is established there (see display ($C.15$)) that
\[
\tau^{-m} \tilde{E}_{m} \left( -\frac{1}{\tau}, \frac{z}{\tau}\right)
=\sum_{k=0}^m \frac{(-1)^{m-k}}{(m-k)!} \tilde{E}_k (\tau ,z) z^{m-k} \tau^{k-m} ,
\]
where we take $\tilde{E}_0 (\tau ,z)$ to be $1$. Therefore,
\begin{equation}
\begin{aligned}
\tau^{-m} \tilde{E}_{m} \left( -\frac{1}{\tau}, \frac{z}{\tau}\right)
&=\sum_{k=0}^m \frac{(-1)^{m-k}}{(m-k)!} \tilde{E}_k (\tau ,z) z^{m-k} \tau^{k-m} \\
&=\sum_{k=0}^m \frac{(-1)^{m-k}}{(m-k)!} \tilde{E}_k (\tau ,z)
\left(\frac{z}{\tau}\right)^{m-k} .
\end{aligned}
\label{quasi-JacobiEisenlemma1}
\end{equation}
 This proves the transformation for the matrix $S$. For the matrix $T$ we have $T\cdot \tau \mapsto \tau +1$, and we find $\tilde{E}_m (\tau +1,z)=\tilde{E}_m (\tau ,z)$. 

 We now consider the general case of $\underline{z}$. Using $\tilde{E}_m (\tau ,\underline{z})=\tilde{E}_m (\tau ,z_1 +\cdots +z_n)$ and (\ref{quasi-JacobiEisenlemma1}), we find
\begin{align*}
  \tau^{-m}\tilde{E}_m  \left(-\frac{1}{\tau}
,\frac{\underline{z}}{\tau}\right)&=
\tau^{-m}\tilde{E}_m  \left(-\frac{1}{\tau} ,\frac{z_1 +\cdots +z_n}{\tau}\right)\\
&=\sum_{k=0}^m \frac{(-1)^{m-k}}{(m-k)!} \tilde{E}_k (\tau ,z_1 +\cdots +z_n) \left(\frac{z_1 +\cdots +z_n}{\tau}\right)^{m-k} \\
&=\sum_{k=0}^m \sum_{\begin{subarray}{c}i_1
,\dots ,i_n \geq 0\\ i_1 +\cdots +i_n =m-k\end{subarray}}  C_{i_1 ,\dots ,i_n}\frac{(-1)^{m-k}}{(m-k)!} \tilde{E}_k (\tau ,\underline{z}) \left(\frac{z_1}{\tau}\right)^{i_1} \cdots \left(\frac{z_n}{\tau}\right)^{i_n},
\end{align*}
where the $C_{i_1 ,\dots ,i_n}$ are scalars produced when expanding $((z_1
+\cdots +z_n)/\tau )^{m-k}$. This proves (\ref{QuasiDefn1}) for the matrix $S$. The case for the matrix $T$ is again trivial.

To prove (\ref{QuasiDefn2}), we can repeat similar steps for $[\lambda ,0] \in \mathbb{Z}^n \times \mathbb{Z}^n$, using the transformation (for the $\underline{z}=z$ case) 
\[
\tilde E_{m}(\tau ,z) =\sum_{k= 0}^m (-1)^{m+k} \binom{m}{k} \lambda^{m-k} \tilde E_{k} (\tau ,z)
\]
(cf.\ \cite{Oberdiek-Serre}, page $7$).\ We omit further details. $\hfill \Box$\\

  The following lemma follows as in Proposition $2$ in \cite{MTZ} with the same proof (see also display (C.14) in \cite{GK-differential}).
 \begin{lemma}\label{BumbleBeeTuna} We have
 \[
 \tilde{P}_m (w,\underline{z},\tau)= \frac{(-1)^m}{w^m}+ \sum_{k\geq m} \binom{k-1}{m-1} \tilde{E}_k (\tau ,\underline{z})w^{k-m}.
 \]
 \end{lemma}
 
 Beyond the modular forms discussed above, we also frequently encounter the usual quasi-modular Eisenstein series $E_2 (\tau)$ normalized so that it has the functional equation
 \[
 E_2 \left(\frac{a\tau +b}{c\tau +d} \right) =(c\tau +d)^2 E_2 (\tau) -\frac{c(c\tau +d)}{2\pi i},
 \]
 for $\gamma =\left(\begin{smallmatrix} a&b \\c&d \end{smallmatrix}\right) \in \text{SL}_2 (\mathbb{Z})$.

\section{Recursion formula\label{Section-RecursionFormula}}

In this section we establish recursion formulas for $n$-point functions. These results are found using
an analysis that is similar to that in \cite{MTZ} and \cite{Zhu}. For this reason, we merely state the needed results, omitting proofs. The next lemma contains the necessary changes as well as the assumption $\zeta_1^{h_1 (0)}\cdots \zeta_m^{h_m (0)}v=v$.
 
\begin{lemma} \label{Sum[0]}
 Let $M^r$ be a module for $V$, $v\in V_{k}$, and $v_1 ,\dots ,v_n \in V$. If $\zeta_1^{h_1 (0)}\cdots \zeta_m^{h_m (0)}v
 \linebreak 
 =v$, then
 \begin{equation}
 \sum_{r=1}^n J_{r,\underline{h}}(v_1 ,\dots ,v[0]v_r ,\dots ,v_n;\tau ,\underline{z})=0.\notag %
 \end{equation}
\end{lemma}

The following two lemmas, when combined, reduce any $n$-point function to a linear combination of 
$(n-1)$-point functions with modular coefficients of the type described in Section \ref{Section-AutomorphicForms}.

\begin{lemma}
 Let $v\in V$ and suppose $h_j (0)v=\mu_j v$, $\mu_j \in \mathbb{C}$, for each $1\leq j\leq m$. Then for any $V$-module $M^r$ and $v_1 ,\dots ,v_n \in V$, we have
 \begin{align}
 &J_{r,\underline{h}}(v,v_1,\dots ,v_n;\tau ,\underline{z}) \notag \\
     &\hspace{5mm}=\delta_{\underline{z}\cdot \underline{\mu},\mathbb{Z}} \tr_{M^r} o(v)Y^M (q_1^{L(0)} v_1,q_1)\cdots Y(q_n^{L(0)}v_n ,q_n) \zeta_1^{h_1 (0)}\cdots \zeta_m^{h_m (0)} q^{L(0)-\mathbf{c}/24}\notag \\
        &\hspace{10mm} +\sum_{s=1}^n \sum_{k\geq 0} \tilde{P}_{k+1} (z_s -z,\tau ,\underline{z}\cdot \underline{\mu}) J_{r,\underline{h}}(v_1 ,\dots ,v[k]v_s ,\dots ,v_n ;\tau ,\underline{z}),
 \notag 
 \end{align}
 where $\delta_{\underline{z}\cdot \underline{\mu},\mathbb{Z}}$ is $1$ if $\underline{z}\cdot \underline{\mu}\in \mathbb{Z}$ and is $0$ otherwise. 
\end{lemma}

\begin{lemma} \label{LemmaRecursion}
 Let the assumptions be the same as in the previous lemma. Then for $p\geq 1$,
 \begin{align}
 &J_{r,\underline{h}}(v[-p]v_1 ,\dots ,v_n;\tau ,\underline{z}) \notag \\
&\hspace{-1mm}= \delta_{\underline{z}\cdot \underline{\mu},\mathbb{Z}}\delta_{p,1}
     \tr_{M^r} o(v)Y^M (q_1^{L(0)} v_1,q_1)\cdots Y(q_n^{L(0)}v_n ,q_n) \zeta_1^{h_1 (0)}\cdots \zeta_m^{h_m (0)} q^{L(0)-\mathbf{c}/24} \notag \\
 &\hspace{0mm} +(-1)^{p+1}\sum_{k\geq 0} \binom{k+p-1}{p-1} \tilde{E}_{k+p}(\tau ,\underline{z}\cdot \underline{\mu}) J_{r,\underline{h}}(v[k]v_1 ,\dots ,v_n ;\tau ,\underline{z})\notag \\
 & +(-1)^{p+1} \sum_{s=2}^n \sum_{k\geq 0} \binom{k+p-1}{p-1}\Big(\tilde{P}_{k+p}(z_s
     -z_1 ,\tau ,\underline{z}\cdot \underline{\mu}) J_{r,\underline{h}}(v_1 ,\dots ,v[k]v_s ,\dots ,v_n;\tau
     ,\underline{z})\Big).
\notag 
 \end{align}
 In particular, in the case $n=1$ we have
  \begin{align}
 J_{r,\underline{h}}(v[-p]v_1;\tau ,\underline{z}) &=\delta_{\underline{z}\cdot \underline{\mu},\mathbb{Z}}\delta_{p,1} \tr_{M^r} o(v)Y^M (q_1^{L(0)} v_1,q_1)\zeta_1^{h_1 (0)}\cdots \zeta_m^{h_m (0)} q^{L(0)-\mathbf{c}/24} \notag \\
 &\hspace{5mm} +(-1)^{p+1}\sum_{k\geq 0} \binom{k+p-1}{p-1} \tilde{E}_{k+p}(\tau ,\underline{z}\cdot \underline{\mu}) J_{r,\underline{h}}(v[k]v_1 ;\tau ,\underline{z}). \notag
\notag 
 \end{align}
\end{lemma}

 Another result that will be useful is the following. (See also \cite{GK-differential}.)
\begin{corollary}
 Let $m=1$ and $n=1$ as in the previous lemma.
 \begin{enumerate}
 \item If $u\in V$ such that $h(0)u=0$, then
 \begin{align}
 J_{r,h}(u[-p]v;\tau ,z) &= \delta_{p,1} \tr_{M^r} o(u)o(v) \zeta^{h (0)} q^{L(0)-\mathbf{c}/24}\notag \\
 &\hspace{5mm} +(-1)^{p+1}\sum_{k\geq 0} \binom{k+p-1}{p-1} E_{k+p}(\tau) J_{r,h}(u[k]v ;\tau ,z). \notag
\notag 
 \end{align}
 \item If $u\in V$ such that $h(0)u=\mu u$ ($\mu \not =0$), then
 \begin{equation}
 J_{r,h}(u[-p]v;\tau ,z)= (-1)^{p+1}\sum_{k\geq 0} \binom{k+p-1}{p-1} \tilde{E}_{k+p}(\tau ,\mu z) J_{r,h}(u[k]v ;\tau ,z). \notag
\notag 
 \end{equation}
 \end{enumerate}
\end{corollary}

\noindent {\bf Remark.}\ \  
 The difference of a minus sign between these equations and those found in \cite{MTZ} can be attributed to the minus sign difference in our definitions of the functions $P_k \left[\begin{smallmatrix} \zeta \\ 1 \end{smallmatrix}\right] (w,\tau)$ and the action of $\text{SL}_2 (\mathbb{Z})$.

 \medskip
 Finally, using that $E_k (\tau)=0$ for odd $k$ we can establish the following corollary.
 
 \begin{corollary} 
  Let $m=n=1$ as before.
 \begin{enumerate} 
 \item For any $v\in V$, we have
     \begin{align}\label{NT1}
         J_{r,h}(h[-1]v;\tau ,z)&=\frac{1}{2\pi i}\frac{d}{dz}J_{r,h}(v;\tau ,z) + \sum_{k\geq 1} E_{2k}(\tau) J_{r,h}(h[2k-1]v ;\tau ,z),
 \end{align}
 and
 \begin{align}\label{NT2}
     J_{r,h}(L[-2]v;\tau ,z)&=\frac{1}{2\pi i}\frac{d}{d\tau}J_{r,h}(v;\tau ,z) + \sum_{k\geq 1} E_{2k}(\tau) J_{r,h}(L[2k-2]v ;\tau ,z).
 \end{align}
 \item If $u\in V$ such that $h(0)u=0$, then
 \begin{equation}
 J_{r,h}(u[-2]v;\tau ,z) = - \sum_{k\geq 1} (2k-1) E_{2k}(\tau) J_{r,h}(u[2k-1]v ;\tau ,z). \label{NT3}
 \end{equation}
 \item If $u\in V$ such that $h(0)u=\mu u$ ($\mu \not =0$), then
 \begin{equation}
 J_{r,h}(u[-2]v;\tau ,z)= \sum_{k\geq 2} (-1)^{k+1} (k-1) \tilde{E}_{k}(\tau ,\mu z) J_{r,h}(u[k-2]v ;\tau ,z). \label{NT4}
 \end{equation}
 \end{enumerate}
\end{corollary}

\section{Theorem \ref{TheoremAlpha} proof\label{Section-TheoremAlpha}}

 Throughout this section, $V$ is a strongly regular vertex operator algebra of central charge $\mathbf{c}$ and $M^1 ,\dots ,M^s$ its inequivalent irreducible admissible modules.\ Fix $h_1 ,\dots ,h_m$ in $V_1$ which satisfy Condition H on each module $M^r$.\ Let $G$ be the Gram matrix $G=(\gby{h_i ,h_j})$ associated with the bilinear form $\gby{\cdot ,\cdot}$ and elements $h_1 ,\dots ,h_m$. 

 We first prove the transformation law (\ref{MainTheorem1}) in Theorem \ref{TheoremAlpha}. To do so we will need a $1$-point analogue of a result due to Miyamoto \cite{Miy}.\ For $u,w\in V_1$ and $v\in V$, we define the function $\Phi_r (v;u,w,\tau)$ by
 \begin{equation}
 \Phi_r (v; u,w,\tau):= \tr_{M^r} o(v) e^{2\pi i (w (0) +\gby{u,w}/2)} q^{L(0)+u(0) +\gby{u,u}/2 -\mathbf{c}/24}.
 \label{phifunction}
 \end{equation}
 Function (\ref{phifunction}) is similar to the functions $\Phi_r$ defined in \cite{Miy}, except there only the case $v=\textbf{1}$ is considered.\ Moreover, we have switched the notation of $u$ and $w$ and taken $\gby{\textbf{1} ,\textbf{1}}=-1$, which is negative the normalization taken by Miyamoto.
 
The proof of the following theorem is the same as in \cite{Miy} (see Theorem A) when one makes the appropriate changes.\ We omit details here; they may be found in \cite{KrauelThesis}.
\begin{theorem}\label{MiyamotoTheorem}
 Let $V$ be a rational, $C_2$-cofinite vertex operator algebra and $M^1 ,\dots ,M^s$ be its finitely many inequivalent irreducible admissible modules. Suppose $w \in V_1$ and $v \in V_{[k]}$ are such that $w(n)v=0$ for $n\geq 0$. Then for all $\gamma =\left( \begin{smallmatrix} a&b\\ c&d \end{smallmatrix} \right) \in \text{SL}_2 (\mathbb{Z})$,
\begin{equation}
\Phi_{r} \left(v;0,w,\frac{a\tau +b}{c\tau +d} \right) =(c\tau +d)^k \sum_{i=1}^s A_{r,\gamma}^i \Phi_{i} (v;cw,dw,\tau),
\notag 
\end{equation}
where $A_{r,\gamma}^i$ are the scalars $S(\gamma ,r,i)$ dependent on $\gamma$ that appear in Zhu's Theorem $5.3.2$ of \cite{Zhu}. 
\end{theorem}
 
Note that
 \[
 J_{r,\underline{h}}(v;\tau ,\underline{z})=\Phi_r (v; 0,\underline{z}\cdot \underline{h},\tau),
 \]
 where $\underline{z}\cdot \underline{h}$ is the usual dot-product. By Theorem \ref{MiyamotoTheorem},
 \begin{align}
 J_{r,\underline{h}}\left(v;\frac{a \tau +b}{c\tau +d} ,\frac{\underline{z}}{c\tau +d} \right)&= \Phi_r \left(v;0, \frac{\underline{z}\cdot \underline{h}}{c\tau +d}, \frac{a \tau +b}{c\tau +d} \right) \nonumber \\
 &= (c\tau +d)^k\sum_{\ell =1}^s A_{r  ,\gamma}^\ell \Phi_\ell \left(v; \frac{c\underline{z}\cdot \underline{h}}{c\tau +d} , \frac{d\underline{z}\cdot \underline{h}}{c\tau +d},\tau \right).
 \label{blergy1}
 \end{align}
 Expanding the $\Phi_\ell$ on the right hand side we find
 \begin{align}
   \Phi_\ell &\left(v; \frac{c\underline{z}\cdot \underline{h}}{c\tau +d} ,
 \frac{d\underline{z}\cdot \underline{h}}{c\tau +d},\tau \right) \notag \\
 &=\tr_{M^\ell} o(v) \exp \left[2\pi i \left(d\frac{\underline{z}\cdot \underline{h (0)}}{c\tau +d} +cd \sum_{j=1}^m \sum_{t=1}^m \frac{ z_j \gby{h_j ,h_t} z_t}{2(c\tau +d)^2} \right)\right] \nonumber \\
  &\hspace{5mm}\cdot \exp \left[ 2\pi i \tau \left(L(0) +c\frac{\underline{z}\cdot \underline{h (0)}}{c\tau +d} +c^2 \sum_{j=1}^m \sum_{t=1}^m \frac{z_j \gby{h_j ,h_t}z_t}{2(c\tau +d)^2} -\mathbf{c}/24\right)\right] \nonumber \\
     &=\tr_{M^\ell} o(v) \exp \left(2\pi i\underline{z}\cdot \underline{h (0)}\right) \exp \left(2\pi i c\sum_{j,t=1}^m \frac{ z_j \gby{h_j ,h_t}z_t}{2(c\tau +d)} \right) q^{L(0)-\mathbf{c}/24}\nonumber \\
     &=\exp \left( \pi i \frac{c G[\underline{z}]}{c\tau +d}\right)
     \tr_{M^\ell} o(v) \zeta_1^{h_1 (0)}\cdots \zeta_m^{h_m (0)}
     q^{L(0)-\mathbf{c}/24} . \label{blergy2}
 \end{align}
 Combining (\ref{blergy1}) and (\ref{blergy2}) establishes (\ref{MainTheorem1}).
 
\medskip
\noindent {\bf Remark.}\ \  
 Although it may appear that Condition H and the assumption $h_j (n)v=0$ are not needed to establish (\ref{MainTheorem1}), they are used in the proof of Theorem \ref{MiyamotoTheorem} and are indeed necessary.

 \medskip
 Next we prove  (\ref{MainTheorem2}).  Following H. Li \cite{Li}, define \emph{invertible} maps 
 $\Delta_{h_j} (z) : V \to (\End V)\left[ z^{-1},z\right]$ by
 \begin{equation}
 \Delta_{h_j} (z) := z^{h_j (0)} \exp \set{-\sum_{k\geq 1} \frac{h_j (k)}{k} (-z)^{-k}}, \notag
 \end{equation}
 and $Y^M_{\Delta_{h_j} (z)} (\cdot ,z) \colon V \to (\End M)\left[
 z^{-1},z\right]$ by
 \[
 Y^M_{\Delta_{h_j} (z)} (v ,z):=Y^M \left( \Delta_{h_j}(z) v,z \right).
 \]
We then have the following theorem (loc.\ cit.\ Proposition $5.4$).
 \begin{theorem}\label{LiTheorem}
 Suppose that $g$ is a finite order automorphism of $V$ such that $g(h_j)=h_j$. Let $(M^r,Y^r)$ be a $g$-twisted $V$-module. Then $(M^r, Y^r_{\Delta_{h_j} (z)})$ is a weak $(ge^{2\pi i h_j (0)})$-twisted $V$-module.
 \end{theorem}
 
 Applying this formalism when $g=e^{2\pi i h_j (0)}$ is the identity automorphism ($h_j$ has 
 \emph{integral} eigenvalues), we obtain an isomorphism of $V$-modules
 \begin{equation}
 (M^{r'} ,Y^{r'}_{\Delta_{h_j} (z)}) \cong (M^r ,Y^{r})
 \label{gtwistedisomorphism}
 \end{equation}
 for some $r' \in \set{1, \dots ,s}$.\ For each $h_i\ (1\leq i \leq m)$ we have
 \begin{align}
 \Delta_{h_i} (z) \omega &= \left(z^{h_i(0)} \exp \set{-\sum_{k\geq 1} \frac{h_i
 (k)}{k} (-z)^{-k}}\right) \omega  \nonumber \\[.5ex]
  &= z^{h_i (0)} \left( \omega -h_i(1)\omega (-z)^{-1} + \frac{h_i(1)^2 \omega}{2} (-z)^{-2} \right)\nonumber \\
  &= \omega + h_i z^{-1} +\gby{h_i ,h_i}z^{-2}.
 \notag 
 \end{align}
 Therefore, the modes of $\omega_{\Delta_{h_i}}$ acting  on $(M^{r'}, Y^{r'}_{\Delta_{h_i} (z)})$ are given by
 \begin{align*}
 \sum_{n\in \mathbb{Z}} \omega_{\Delta_{h_i}} (n) z^{-n-1} &=Y^{r'}_{\Delta_{h_i} (z)} (\omega ,z)=Y^{r'}(\Delta_{h_i}(z)\omega ,z)\\
  &=\sum_{n\in \mathbb{Z}} \left( \omega (n)z^{-n-1} + h_i(n) z^{-n-2} +\frac{\gby{h_i ,h_i}}{2} z^{-n-3} \right).
 \end{align*}
 Taking $\text{Res}_{z} z$ of both sides, we find $\omega_{\Delta_{h_i}}(1) =\omega (1) +h_i (0) +\gby{h_i ,h_i}/2$, i.e.,
 \begin{equation}
 L_{\Delta_{h_i}} (0)= L(0)+h_i (0) +\frac{\gby{h_i ,h_i}}{2}.
 \label{blergith2}
 \end{equation}
 In a similar way we have
 \[
 \Delta_{h_i}(z) h_j = h_j +\gby{h_i ,h_j}z^{-1}
 \]
 for any $1\leq j \leq m$, and in particular,
 \begin{equation}
 (h_j)_{\Delta_{h_i}}(0)= h_j (0) +\gby{h_i ,h_j}.
 \label{blergith3}
 \end{equation}
Using (\ref{blergith2}) and (\ref{blergith3}), we  find
 \begin{equation}
 L_{\Delta_{-\underline{\lambda}\cdot \underline{h}}}(0)=L(0) -\underline{\lambda}\cdot \underline{h} +\frac{1}{2}\sum_{s=1}^m \sum_{t=1}^m \lambda_s \gby{h_s ,h_t}\lambda_t , \label{blergith4} 
 \end{equation}
 and
 \begin{equation}
 (h_j)_{\Delta_{-\underline{\lambda}\cdot \underline{h}}} (0)= h_j (0) -\sum_{s=1}^m \lambda_s \gby{h_s ,h_j}. \label{blergith5}
 \end{equation}
 Finally, with these same calculations applied to $v\in V_{[k]}$ such that $h_j (n)v=0$ for $1\leq j \leq m$ and $n\geq 0$, we find that the zero mode of $v$ on $(M^{r'} ,Y^{r'}_{\Delta_{-\underline{\lambda}\cdot \underline{h}} (z)})$ is $o_{\Delta_{-\underline{\lambda}\cdot \underline{h}}} (v) = v(k-1)=o(v)$.
 
 Using (\ref{blergith4}), (\ref{blergith5}), and the isomorphism (\ref{gtwistedisomorphism}), it follows that
 \begin{align*}
J_{r,\underline{h}} &(v;\tau ,\underline{z}+\underline{\lambda}\tau +\underline{\mu})\\
&=\tr_{M^r} o(v) \exp \left[2\pi i \sum_{\delta=1}^m \left(z_\delta
 +\lambda_\delta \tau +\mu_\delta \right)h_\delta (0) \right] \exp (2\pi i \tau
 (L(0)-\mathbf{c}/24))\\
  &=\tr_{M^r} o(v) \exp \left[2\pi i \sum_{\delta=1}^m \left(z_\delta +\lambda_\delta \tau  \right)h_\delta (0) \right] \exp (2\pi i \tau (L(0)-\mathbf{c}/24))\\
 &=\tr_{M^{r'}} o(v) \exp \left[ 2\pi i \sum_{\delta=1}^m \left(z_\delta +\lambda_\delta \tau\right) \left(h_\delta (0) -\sum_{t=1}^m \lambda_t \gby{h_t ,h_\delta} \right)\right]\\ 
 &\quad\ \cdot \exp \left[ 2\pi i \tau \left(L(0)-\underline{\lambda}\cdot \underline{h(0)} +\frac{1}{2}\sum_{\delta=1}^m \sum_{t=1}^m \lambda_\delta \gby{h_\delta ,h_t} \lambda_t -\mathbf{c}/24 \right) \right]\\
 &=\tr_{M^{r'}} o(v) \exp \left(2\pi i \underline{z}\cdot \underline{h(0)}\right) \exp \left(2\pi i \tau \underline{\lambda}\cdot \underline{h(0)} \right) \\
 &\quad\ \cdot \exp \left(-2\pi i \sum_{\delta=1}^m \sum_{t=1}^n z_\delta \gby{h_\delta ,h_t}\lambda_t \right)  \exp \left(-2\pi i \tau \sum_{\delta=1}^m \sum_{t=1}^m \lambda_\delta \gby{h_\delta ,h_t} \lambda_t \right) \\
 &\quad\ \cdot \exp \left( -2\pi i \tau \underline{\lambda}\cdot \underline{h(0)}\right)  \exp \left(\pi i \tau \sum_{\delta=1}^m \sum_{t=1}^m \lambda_\delta \gby{h_\delta ,h_t}\lambda_t \right) q^{L(0)-\mathbf{c}/24}\\
 &=\exp \left(-\pi i \left(G[\underline{\lambda}]\tau + 2\underline{z}^t G\underline{\lambda} \right)\right) J_{r',\underline{h}}(v;\tau ,\underline{z}).
 \end{align*}
Here,  the second equality uses the fact $\exp \big(2\pi i \underline{\mu}\cdot
\underline{h(0)} \big) =1$ since $\underline{\mu}\cdot \underline{h(0)}$ acts on $M^r$ with integer eigenvalues.\ This proves (\ref{MainTheorem2}).

 Since (\ref{FourierExpansion}) is clear, it remains to establish the convergence of $J_{r, \underline{h}}(v;\tau ,\underline{z})$ for any $v\in V$ and module $M^r =\bigoplus_{d\geq 0} M_{\lambda_r +d}^r$, where $\lambda_r$ is the conformal weight of $M^r$. For the remainder of this section we also drop the notation $J_{r,h}$ and simply write $J_{h}$ as none of the calculations are dependent on $M^r$.

 Consider the case $m=1$. That is, take $\underline{h}$ to be a single element $h\in V_1$ that satisfies Condition H, so that we are concerned with the function $J_{r,h}(v;\tau ,z)$ on $\mathbb{H}\times \mathbb{C}$. Let $\mathcal{M}$ denote the ring of quasi-modular forms and $\mathcal{Q}_0$ be the ring of quasi-Jacobi forms of index $0$ (see \cite{Libgober-Elliptic}, Proposition $2.8$), which are both known to be Noetherian. By Definition $2.5$ in \cite{Libgober-Elliptic}, it is clear certain generators of quasi-Jacobi forms of index $0$ are convergent on closed subsets of $\{ (\tau ,z) \in \mathbb{H}\times \mathbb{C} \mid z\not \in \mathbb{Z}+\mathbb{Z}\tau \}$. In particular, the functions $\tilde{E}_k(\tau ,z)$ introduced in Subsection \ref{Subsection-TwistedEllipticFunctions} are convergent on this domain.

 \enlargethispage{1em}
 For $\alpha \in \mathbb{Z}$, let $U_\alpha$ be the map defined by $U(\phi (\tau ,z))=\phi (\tau ,\alpha z)$. Then $U_\alpha$ maps (quasi-)Jacobi forms of index $m$ to (quasi-)Jacobi forms of index $\alpha^2 m$. In particular, $\tilde{E}_{k} (\tau ,\alpha z) =U_\alpha (\tilde{E}_{k}(\tau ,z))$ is a quasi-Jacobi form of index $0$.

 Set $V(\mathcal{Q}_0) =V\otimes \mathcal{Q}_0$ and let $O_h(V)$ be the subspace of $V$ generated by the elements
 \begin{align}
     &u[0]v, \\
     &u[-2]v +\sum_{k=2}^\infty (2k-1) E_{2k}(\tau)u[2k-1]v,\quad \text{when }h(0)u=0, \text{  and}\\
     &u[-2]v +\sum_{k=2}^\infty (k-1)\tilde{E}_{k}(\tau ,\alpha
     z)u[k-2]v, \quad\text{when }h(0)u=\alpha z.
 \end{align}
 By (\ref{NT3}), (\ref{NT4}), and Lemma (\ref{Sum[0]}), it follows that $J_h (v,\tau ,z)=0$ for all $v\in O_h(V)$.
 
 \begin{lemma}
 Suppose $V$ is $C_2$-cofinite. Then $V(\mathcal{Q}_0)/O_h(V)$ is a finitely generated $\mathcal{Q}_0$-module.
 \end{lemma}
 
 \paragraph{Proof}
 The proof mimics that of Lemma $4.4.1$ in \cite{Zhu}. Since $C_2(V)$ has finite codimension, there exists an integer $N$ such that $V_n \subset C_2 (V)$ for all $n>N$. Let $A$ be the $R$-submodule of $V(\mathcal{Q}_0)$ generated by $\bigoplus_{n\leq N} V_n$. If $v\in V_{[k]}$ we will show that $v\in A+O_h(V)$, thereby proving $V(\mathcal{Q}_0)=A+O_h(V)$, and thus the lemma. 

 In the case $k\leq N$, we are done since $v\in A$. Therefore we assume that $k>N$. In this case, $V_{[k]} \subset A +O_h(V)$, and so we have $v =a +\sum_{i=0}^d b_i(-2)c_i$ for some $a\in A$ and homogeneous $b_i,c_i \in V$ satisfying $\wt [b_i]+\wt [c_i] =k -1$. In the case $h(0)b_i =0$, then $b_i (-2)c_i \in A +O_h(V)$ just as in \cite{Zhu}. It suffice to show $b_i (-2)c_i \in A +O_h(V)$ in the case $h(0)b_i =\alpha b_i$ for some nonzero $\alpha \in \mathbb{Z}$.

 In this case we have
 \[
 b_i [-2] c_i + \sum_{\ell =2}^\infty (\ell -1) \tilde{E}_{\ell}(\tau ,\alpha z) b_i [\ell -2] c_i
 \]
 is in $O_h(V)$. Since $\tilde{E}_{k}(\tau ,\alpha z)$ is again a quasi-Jacobi form and $\wt [b_i [\ell-1]c_i] =\wt [b_i]+\wt [c_i] -\ell =k-1-\ell$, our induction hypothesis shows that
 \[
 \sum_{\ell =2}^\infty (\ell -1) \tilde{E}_{\ell}(\tau ,\alpha z)b_i [\ell -2]c_i
 \]
 is in $O_h(V)$, and thus so is $b_i [-2] c_i$. Using that $b_i (-2)c_i =b_i[-2] c_i + \sum_{j>-2} \beta_j b_i [j]c_i$ for some scalars $\beta_j$, we can apply our induction hypothesis again to the elements $b_i [j] c_i$ to find $b_i (-2) c_i \in O_h (V)$. The lemma is now proved. $\hfill \Box$\\
 
\noindent {\bf Remark.}\ \   
Since our recursion formula introduces functions of the form $\tilde{E}_\ell (\tau ,\alpha z)$ ($\alpha \in \mathbb{Z}$) in the previous step, an arbitrary function $\phi (\tau ,z)$ in $\mathcal{Q}_0$ may have a pole at $z \in \frac{1}{\alpha}\left( \mathbb{Z} +\mathbb{Z}\tau \right)$ for different $\alpha$. Therefore, the finite many elements in $\mathcal{Q}_0$ that arise in the following lemma may each have such poles. The authors would like to thank Reimundo Heluani for bringing this to their attention. After this step, however, no further coefficients with poles are introduced in the proof, and we therefore obtain the domain of convergence as described in the statement of the theorem.
 
 \begin{lemma}
 Suppose $V$ is $C_2$-cofinite. For any $v\in V$ there exist $m,n\in \mathbb{N}$ and $\phi_i (\tau ,z), \psi_j (\tau ,z) \in \mathcal{Q}_0$, $0\leq i\leq m-1$, $0\leq j\leq n-1$, such that
 \begin{equation}
 L[-2]^m v +\sum_{i=0}^{m-1} \phi_i (\tau ,z) L[-2]^i v \in O_h (V) \label{TryNewEq1}
 \end{equation}
 and
 \begin{equation}
 h[-1]^n v +\sum_{j=0}^{n-1} \psi_j (\tau ,z) h[-1]^j v \in O_h (V). \label{TryNewEq2}
 \end{equation}
 \end{lemma}
 
 \paragraph{Proof}
 By the previous lemma and the fact $\mathcal{Q}_0$ is Noetherian, we have the $\mathcal{Q}_0$-submodule generated by $\{h[-1]^j v, j\geq 0 \}$ is finitely generated. Therefore, some relation such as (\ref{TryNewEq2}) must hold. Equation (\ref{TryNewEq1}) is proved similarly. $\hfill \Box$\\
 
 Set $D_\tau =\frac{1}{2\pi i}\frac{d}{d\tau}$ and $D_z =\frac{1}{2\pi i}\frac{d}{dz}$.
 
\begin{proposition}\label{NewTryProp}
 Suppose that $V$ is $C_2$-cofinite.
 \begin{enumerate}
 \item If $v\in V$ is such that $L[\ell]v =0$ for $\ell >0$, then there exists an $m\in \mathbb{N}$ and $\phi_i (\tau ,z) \in \mathcal{Q}_0$, $0\leq i \leq m -1$, such that
 \begin{equation}
 D_\tau^m J_h (v,\tau ,z) +\sum_{i=0}^{m-1} \phi_i (\tau ,z) D_\tau^i J_h (v,\tau ,z) =0. \label{NewTry3}
 \end{equation}
 \item If $v\in V$ is such that $h[\ell]v =0$ for $\ell >0$, then there exists an $n\in \mathbb{N}$ and $\psi_j (\tau ,z) \in \mathcal{Q}_0$, $0\leq j \leq n -1$ such that
 \begin{equation}
 D_z^n J_h (v,\tau ,z) +\sum_{j=0}^{n-1} \psi_j (\tau ,z) D_z^j J_h (v,\tau ,z) =0. \label{NewTry4}
 \end{equation}
 \end{enumerate}
\end{proposition} 

\paragraph{Proof}
 The proof of (\ref{NewTry3}) follows just as in \cite{DLM-orbifold}. The proof of (\ref{NewTry4}) is similar, and follows from using Equation (\ref{NT1}) along with induction and the fact that there are scalars $\beta_{ijk}$ such that $h[2k-1]h[-1]^i v=\sum_{j=0}^{i-1} \beta_{ijk} h[-1]^j v$ for any $k\geq 1$. $\hfill \Box$\\
 
 It follows from the theory of ordinary differential equations that solutions to these equations converge wherever the functions $\phi_i (\tau ,z)$ and $\psi_j (\tau ,z)$ do. In particular, (\ref{NewTry3}) shows that $J_h (v,\tau ,z)$ converges on the set $F:=\{ (\tau ,z)\in \mathbb{H}\times \mathbb{C} \mid z \not \in \frac{1}{t}(\mathbb{Z}+ \mathbb{Z}\tau), \text{ finitely many }t\in \mathbb{Z} \}$.

It remains to show that similar differential equations hold for any $v\in V$,
not just those that are primary. We will first establish a series of lemmas
involving a single variable $z$. We omit the variable $\tau$ until the end, as
the analogous results can be proved similarly and are also essentially found in~\cite{DLM-orbifold}.
 
\begin{lemma}
 Suppose $\ell \geq 1$, $j\geq 0$, and $v\in V_{[k]}$ is such that $h(0)v=\alpha v$ for some scalar $\alpha$.
 \begin{enumerate}
 \item For $\ell =1$, $h[\ell -1]h[-1]^j v =\alpha h[-1]^j v$.
 \item For $\ell =2$, $h[\ell -1]h[-1]^j v = h[-1]^{j} h[1]v + j\langle h,h\rangle h[-1]^{j-1} v$.
 \item For all $\ell \geq 1$, there are scalars $\beta_{ij\ell}$ and elements satisfying $\wt[u_{ij\ell}]\leq \wt [v]$, with equality only if $u_{ij\ell} =v$, such that
 \[
 h[\ell -1]h[-1]^j v = h[-1]^j h[\ell -1]v + \sum_{i=0}^{j-1} \beta_{ij\ell} h[-1]^i u_{ij\ell}.
 \]
 \end{enumerate}
\end{lemma}

 \paragraph{Proof}
 Both (a) and (b) follow from easy proofs by induction on $j$. Part (c) follows by induction on $j+\ell$. $\hfill \Box$\\
 
 Using this lemma along with (\ref{NT1}) we find
 \begin{align}
 J_h (h[-1]^{j+1} v,\tau ,z) &= D_z J_h (h[-1]^{j} v,\tau ,z) +\sum_{\ell =1}^\infty E_{2\ell}(\tau) J_h (h[2\ell -1]h[-1]^{j}v,\tau ,z)  \notag \\
     &= D_z J_h (h[-1]^{j} v,\tau ,z)
 +\sum_{\ell =1}^\infty E_{2\ell}(\tau) \Bigg( J_h (h[-1]^{j} h[2\ell -1]v,\tau
 ,z) \notag\\
 &\hspace{5mm}+\sum_{i =0}^j \beta_{ij\ell} J_h (h[-1]^{i} u_{ij\ell},\tau
 ,z)\Bigg).\label{NewTry5}
 \end{align}
 Noting that $\wt [h[2\ell -1]v] < \wt [v]=k$ for $\ell \geq 1$ and using induction on $k+j$ along with (\ref{NewTry5}), we obtain the following lemma. 
 
  \begin{lemma}
 For any $v\in V_{[k]}$ and $j\geq 0$ there are elements $u_{ij\ell} \in V$ satisfying $\wt [u_{ij\ell}] < k$ and functions $f_{ij}(\tau), g_{ij\ell}(\tau) \in \mathcal{M} \subset \mathcal{Q}_0$, $0\leq i\leq j-1$, such that
 \begin{equation}
 \begin{aligned}
 J_{h}(h[-1]^j v,\tau ,z)
 &= D_z^j J_{h} (v,\tau ,z) +\sum_{i=0}^{j-1} f_{ij}(\tau) D_z^i J_h (v,\tau
 ,z)  \\
 &\quad +\sum_{i=0}^{j-1} \sum_{\ell} g_{ij\ell} (\tau) D_z^i J_h (u_{ij\ell},\tau ,z). 
 \end{aligned}\label{NewTry6}
 \end{equation}
 \end{lemma}
 
 \begin{proposition}
 For any $v\in V_{[k]}$ there exist $m,n \in \mathbb{N}$, elements $u_{i\ell},
 w_{j\ell} \in V$ satisfying $\wt [u_{i\ell}],\wt [u_{i\ell}']  < k$, and
 functions $\phi_{i}(\tau ,z)$, $\phi_{i\ell}(\tau ,z)$, $\psi_{j}(\tau ,z)$, $\psi_{j\ell}(\tau ,z) \in \mathcal{Q}_0$, $0\leq i\leq m-1$, $0\leq j\leq n-1$, such that
 \begin{align}\label{NewTryLemmaE2}
     D_\tau^m J_{h} (v,\tau ,z) &+\sum_{i=0}^{m-1} \phi_{i}(\tau ,z) D_\tau^i J_h
 (v,\tau ,z)+\sum_{i=0}^{m-1} \sum_{\ell} \phi_{i\ell} (\tau ,z) D_\tau^i J_h
 (u_{i\ell},\tau ,z) =0 
 \end{align}
 and
 \begin{align}\label{NewTryLemmaE1}
     D_z^n J_{h} (v,\tau ,z) &+\sum_{j=0}^{n-1} \psi_{j}(\tau ,z) D_z^{j} J_h
     (v,\tau ,z)+\sum_{j=0}^{n-1} \sum_{\ell} \psi_{j\ell} (\tau ,z) D_z^{j} J_h
     (w_{j\ell},\tau ,z) =0. 
     \end{align}
 (Note that the functions $\phi$ and $\psi$ are not necessarily the same as
 those in~(\ref{TryNewEq2}).)
 \end{proposition}
 
 \paragraph{Proof}
 Since $J_h (u,\tau ,z)=0$ for $\forall u\in O_h(V)$, we can substitute (\ref{TryNewEq2}) into $J_h$ to obtain $0$. Next, solve for $J_h (h[-1]^n v,\tau ,z)$ and then exchange $J_h (h[-1]^n v,\tau ,z)$ with the right hand side of (\ref{NewTry6}), replacing $j$ with $n$. What results is (\ref{NewTryLemmaE1}). Equation (\ref{NewTryLemmaE2}) follows from a similar analysis, where the necessary lemmas analogous to those here can be proved just as in \cite{DLM-orbifold}. $\hfill \Box$\\
 
 We are now in position to prove the convergence of $J_h (v,\tau ,z)$ for any $v\in V$.\ We first fix $\tau$ and prove that $J_h (v,\tau ,z)$ converges in the $z$-variable on the set $F$.\ The same proof can be applied to prove that $J_h (v,\tau ,z)$ converges for all $\tau \in \mathbb{H}$ by fixing $z$, though we will omit these details.\ We proceed by induction on $\wt [v]$ for homogeneous elements $v\in V$.\ If $\wt [v]=0$, then $v=\beta \textbf{1}$ for some scalar $\beta$, and $J_h (v,\tau ,z)$ satisfies the relevant differential equation in Proposition (\ref{NewTryProp}), and therefore converges on $F$ since the functions $\psi_j (\tau ,z)$ do.\ Suppose, then, that for any $v\in V$ with $\wt [v] \leq k-1$, $J_h (v,\tau ,z)$ also converges on this domain, and consider the case $\wt [v]=k$.

 By our induction hypothesis, the functions $J_h (w_{j\ell},\tau ,z)$ in the previous proposition all converge on $F$.\ For the same fixed $\tau$, set $G(\tau ,z)$ to be the third summand in (\ref{NewTryLemmaE1}). That is, $G(\tau ,z)=\sum_{j=0}^{n-1} \sum_{\ell} \psi_{j\ell} (\tau ,z) D_z^j J_h (w_{j\ell},\tau ,z)$.\ Then a well-known result (see for example Lemma $1$ in \cite{AM}) asserts the existence of a function $k(z)$ that converges on the same domain $F$ and satisfies $(D_z +k(z))G(\tau ,z)=0$.\ Applying the operator $D_z +k(z)$ to (\ref{NewTryLemmaE1}) shows that $J_h (v,\tau ,z)$ satisfies a differential equation with respect to the operator $D_z$ and with coefficients that converge on $F$ for fixed $\tau$.\ This in turn implies $J_h (v,\tau ,z)$ converges on $F$.\ As mentioned before, fixing $z$ and using the same argument with (\ref{NewTryLemmaE2}) shows $J_h (v,\tau ,z)$ also converges for all $\tau \in \mathbb{H}$.\ This proves $J_h (v,\tau ,z)$ converges on the domain stated in Theorem (\ref{TheoremAlpha}) in the case $m=1$.
 
To prove the convergence for the function $J_{r,\underline{h}}(v;\tau ,\underline{z})$ when $m>1$, we fix all but one of the complex variables $z_1 ,\dots ,z_m$ and apply the previous argument.\ Since the convergence can be established in this manner for each individual complex variable, Hartog's Theorem  gives the convergence of $J_{r,\underline{h}}(v;\tau ,\underline{z})$.\ The proof of Theorem \ref{TheoremAlpha} is complete.

\section{Theorem \ref{TheoremBeta} proof\label{Section-TheoremBeta}}

In this section we take up the transformation laws of the functions $J_{r,\underline{h}}(v;\tau ,\underline{z})$ when $h_j (n)v\not =0$ for some $1\leq j \leq m$ or $n\geq 0$. We begin by reviewing the Heisenberg VOA and a decomposition for strongly rational VOAs.

\subsection{The Heisenberg VOA and a module decomposition}

Let $H$ be a $d$-dimensional abelian Lie algebra with non-degenerate symmetric
invariant bilinear form $(\cdot ,\cdot)$. Consider the affinization
$\widehat{H}=H\otimes \mathbb{C}[t, t^{-1}] \oplus \mathbb{C}K$, where $K$ is  central and $[a \otimes t^m, b \otimes t^n]=(a,b) \delta_{m+n,0} K\
(a,b \in H, m,n\in \mathbb{Z})$. Taking $K$ to act as $1$ on $\mathbb{C}$ and $H\otimes \mathbb{C}[t]$ to act trivially, we produce the induced module $M_H =\mathcal{U}(\widehat{H}) \otimes_{H\otimes \mathbb{C}[t]\oplus \mathbb{C}K} \mathbb{C}$ which is isomorphic to the symmetric algebra $S(H\otimes t^{-1} \mathbb{C}[t^{-1}])$ as linear spaces.

 Let the action of $u\otimes t^n$ on $M_H$ be denoted by $u(n)$. For an orthonormal basis $\set{u_1 ,\dots ,u_d}$ of $H$, set $\omega_{M_H} =\frac{1}{2}\sum_{i=1}^d u_i (-1)^2 \textbf{1}$, where $\textbf{1}=1\otimes 1$. Any element $v\in M_H$ can be written as a linear combination of elements of the form
 \[
 v=a_1 (-n_1)\cdots a_\nu (-n_\nu)\textbf{1},
 \] 
 for $a_1 ,\dots ,a_\nu \in H$ and $n_1, \dots ,n_\nu \in \mathbb{N}$. For such
 an element $v\in M_H$, define the map $Y(\cdot ,z)\colon M_H \to (\End
 M_H)\left[ z,z^{-1}\right]$ by
 \[
 Y(v,z)=\normordprod{\partial^{(n_1 -1)} a_1 (z)\cdots \partial^{(n_\nu -1)} a_\nu (z)},
 \]
 where $\partial^{(n)}=\frac{1}{n!}\left(\frac{d}{dz}\right)^n$, $a_i (z) =\sum_{n\in \mathbb{Z}} a_i (n) z^{-n-1}$ ($1\leq i\leq \nu$), and $\normordprod{\cdots}$ signifies normal ordering (see for example \cite{LL}).

 It is known that $(M_H,Y,\textbf{1},\omega_{M_H})$ is a simple (though not rational) vertex operator algebra of central charge $d$ with $L(0)$-grading
 \[
 M_H =\bigoplus_{n\geq 0} (M_H)_n ,
 \]
 where 
 \begin{align*}
 (M_H)_n &=\set{v\in M_H \mid L(0)v=nv} \\
 &=\Big\langle a_1 (-n_1)\cdots a_\nu (-n_\nu )\textbf{1} \mid a_1 ,\dots ,a_\nu
     \in H, n_1 ,\dots ,n_\nu \in \mathbb{N},\ \sum n_i =n\Big\rangle.
 \end{align*}
 There is a natural identification between $(M_H)_1$ and $H$ given by $u(-1)\textbf{1} \mapsto u$. Moreover, for $a,b\in H$ we have $a[0]=a(0)=0$ and $a[1]b=a(1)b=(a,b)\textbf{1}$.

 For $\alpha \in H$, define the space
 \begin{equation}
 M_H (\alpha):= M_H \otimes e^\alpha . \notag 
 \end{equation}
 If $n\not =0$ the operators $a(n) \in \End M_H$ act on $M_H (\alpha)$ via its action on $M_H$.\ On the other hand, $a(0)$ acts on $e^\alpha$ by $a(0)e^\alpha =(a,\alpha)e^\alpha$.\ The space $M_H (\alpha)$ is an irreducible $M_H$-module with conformal weight $\frac{1}{2}(\alpha ,\alpha)$, and
 for varying $\alpha$ we obtain in this way all of the irreducible $M_H$-modules up to equivalence (see \cite{LL} for details).

 The partition function $Z_{M_H}(\textbf{1},\tau):=\tr_{M_H}q^{L(0)-d/24}$ for $M_H$ satisfies
 \begin{equation}
 Z_{M_H}(\textbf{1},\tau)=\eta (\tau)^{-d}.
\notag 
 \end{equation}
 Therefore, since $L(0)e^\alpha =\frac{1}{2}(\alpha ,\alpha)e^\alpha$, $h_j (0)e^{\alpha}=(h_j ,\alpha) e^\alpha$, and $h_j (0)M_H =0$ for all $1\leq j\leq m$, we find
 \begin{align}
    &J_{M_H (\alpha),\underline{h}}(\textbf{1};\tau ,\underline{z})=\tr_{M_H \otimes e^\alpha} \zeta_1^{h_1 (0)}\cdots \zeta_m^{h_m (0)}q^{L(0)-d/24} \notag \\
 &\hspace{5mm}=\left(\tr_{M_H} \zeta_1^{h_1 (0)}\cdots \zeta_m^{h_m (0)}q^{L(0)-d/24}\right) \left(\tr_{e^\alpha} \zeta_1^{h_1 (0)}\cdots \zeta_m^{h_m (0)}q^{L(0)-d/24}\right)\notag \\
 &\hspace{5mm}=\left(Z_{M_H}(\textbf{1},\tau) \right) \zeta_1^{(h_1 ,\alpha)} \cdots \zeta_m^{(h_m,\alpha)} q^{\frac{1}{2}(\alpha ,\alpha)} \notag \\
 &\hspace{5mm}= \eta (\tau)^{-d}\zeta_1^{(h_1 ,\alpha)} \cdots \zeta_m^{(h_m,\alpha)}
     q^{\frac{1}{2}(\alpha ,\alpha)}.\label{HeisenbergJpartition}
 \end{align}
 
 We now discuss a decomposition for any irreducible $V$-module $M$ (see \cite{Mason-Lattice} for more details).\ A result of Dong and Mason \cite{DM-Effective,Mason-Lattice} states that $V_1$ is a \emph{reductive} Lie algebra and that $M$ is a linearly reductive $V_1$-module,
 i.e.,  its action on $M$ is completely reducible. The action of $u \in V_1$ on $M$ is given by $u(0)$.\ Let 
\[
\Omega_M := \set{w\in M \mid u(n)w=0, \text{ for }u\in H \text{ and }n\geq 1},
\]
 and for $\beta \in H$ set
 \[
 M(\beta) := \set{w\in M \mid u(0)w=(\beta ,u)w, \text{ where }u\in H}.
 \]
 Consider the set
 \[
 P:= \set{\beta \in H \mid M(\beta) \not =0},
 \]
 which is a subgroup of $H$. Then $M$ has a decomposition
 \begin{equation}
 M\cong M_H \otimes \Omega_M =\bigoplus_{\beta \in P} M_H \otimes \Omega_M (\beta),
 \label{Decomp1}
 \end{equation}
 where $\Omega_M (\beta):= \Omega_M \cap M(\beta)$ (cf.\ \cite{DM-Integrability,Mason-Lattice}). 

 It is known that $\Omega_V (0)=\Omega (0)$ is a simple vertex operator algebra and $\Omega_M (\beta)$ are irreducible $\Omega (0)$-modules. Moreover, we have $M_H (\beta) \cong M_H \otimes e^\beta$, where $e^\beta \in \Omega (\beta)$. It follows that the tensor product $M_H (\beta) \otimes \Omega_M (\beta)$ is an irreducible $M_H \otimes \Omega (0)$-module. Note also that $M(\beta)=M_H (\beta) \otimes \Omega_M (\beta)$.

 Set
 \[
 L_0 =\set{u\in H \mid u(0) \text{ as an operator on }M \text{ has eigenvalues in }\mathbb{Z}},
 \]
 and 
 \begin{equation}
 \Lambda := \set{u\in L_0 \mid (M,Y^M_{\Delta_u (z)})\cong (M,Y^M)}.
 \notag 
 \end{equation}
 Then the isomorphism (\ref{gtwistedisomorphism}) (which holds for all $u\in L_0$) implies
 \begin{equation}
 \Omega_M (\beta) \cong \Omega_M (\beta +u),
 \label{OmegaIso}
 \end{equation}
 where $u\in \Lambda$ and $\beta \in P$.\ In the case $\beta =0$, this gives $\Omega_M (u) \cong \Omega_M (0)$ for all $u\in \Lambda$.\ 
 Therefore, $\Omega_M (u) \not =0$ and $\Lambda \subseteq P$.\ In \cite{Mason-Lattice}, it is shown that $\Lambda$ is a positive-definite integral lattice of rank $d$ and $|P \colon \Lambda|$ is finite. We set $\delta :=|P \colon \Lambda|$.

 The decomposition (\ref{Decomp1}), which is an isomorphism of vector spaces, may now be written as a decomposition of irreducible modules for $M_H (0)\otimes \Omega_M (0)$. Namely,
 \begin{align}
 M&=\bigoplus_{t=1}^\delta \bigoplus_{\beta \in \Lambda} M_H (\beta +\gamma_t) \otimes \Omega_M (\gamma_t) \notag \\
 &=\bigoplus_{t=1}^\delta \bigoplus_{\beta \in \Lambda +\gamma_t} M_H (\beta) \otimes \Omega_M (\beta), \label{Decomp2}
 \end{align}
 where $\set{\gamma_t }$ are coset representatives of $P/\Lambda$.

 \subsection{Proof of Theorem \ref{TheoremBeta}\label{Subsection-TheoremBeta}}
 
 Since $M_H$ is a subspace of $V$, $H$ is also a subspace of $V$ ( identifying  $a(-1)\textbf{1}\in (M_H)_1$ with $a\in H$). Therefore, we may take the bilinear form $(\cdot ,\cdot)$ on $H$ considered in the previous section to be one which, when extended to $V$, is the restriction of the bilinear form $\gby{\cdot ,\cdot}$ (on $V$) to $H$. In other words, we have $(\cdot ,\cdot)=\gby{\cdot ,\cdot}$ on $M_H$ and we will fix $(\cdot ,\cdot)$ to be such a bilinear form on $H$ and use the notation $\gby{\cdot ,\cdot}$ for the remainder of the paper.

 Let $\set{u_i \mid 1\leq i \leq d}$ be a basis for $H$.\ By the decomposition (\ref{Decomp2}) of $V$, any element in $V$ may be written as sums of elements of the form
 \begin{align*}
     v=u_1 [-m_{1,1}]^{\ell_{1,1}} \cdots u_1 [-m_{1,\nu_1}]^{\ell_{1,\nu_1}}
     &\cdots  u_d [-m_{d,1}]^{\ell_{d,1}} \cdots u_d [-m_{d,\nu_d}]^{\ell_{d,\nu_d}} \otimes e^\alpha \otimes w,
 \end{align*}
  $w\in \Omega (\alpha)$, for various $\alpha \in \Lambda +\gamma_t$, $1\leq
  t\leq \delta$ and $\ell_{x,y} ,m_{x,y} \in \mathbb{N}$ ($1\leq x \leq d$,
  $1\leq y\leq \nu_d$). Note that $v(n) M_H (\beta)\otimes \Omega_{M^r} (\beta)
  \subseteq M_H (\alpha +\beta)\otimes \Omega_{M^r} (\alpha +\beta)$ for an
  irreducible $V$-module $M^r$. Therefore, the only $v$ such that $J_{r,\underline{h}}(v;\tau ,\underline{z}) \not =0$, are those that are a sum containing terms which lie in $M_H (0)\otimes \Omega (0)$.\ That is, for $\alpha =0$ and $w\in \Omega (0)$. It therefore suffices to consider elements of the form given in (\ref{introeq2}).

 Since $w\in \Omega (0)$, it satisfies $h_j (0)w=\gby{h_j ,w}w=0$ for all $1\leq j \leq m$, and $h_j (n)w=0$ for all $n\geq 0$.\ Therefore, $J_{r,\underline{h}}(w;\tau ,\underline{z})$ satisfies the assumptions of Theorem \ref{TheoremAlpha}.

 We will prove Theorem \ref{TheoremBeta} for $v$ as in (\ref{introeq2}) by first establishing results for specific $v$.\ The following lemma and proof follow those found in \cite{DMN-quasi}.
\begin{lemma}\label{quasi1}
 Let $a\in M_H$.\ Consider an element $a [-1]^\ell w \in V$, $\ell \geq 0$, $w\in \Omega (0)$, and let $\alpha \in \Lambda +\gamma_t $ for some $1\leq t \leq \delta$.\ Then there are scalars $c_{\ell ,\ell -2i}$ with $0\leq i \leq \ell /2$ and $c_{\ell ,\ell}=1$ such that
 \begin{align}
 &J_{M_H (\alpha)\otimes \Omega_{M^r} (\alpha), \underline{h}}(a [-1]^\ell w ;\tau ,\underline{z}) \notag \\
 &\hspace{5mm}=\left( \sum_{0\leq i \leq \ell /2} c_{\ell ,\ell -2i} \gby{a,\alpha}^{\ell -2i} (\gby{a,a}E_2 (\tau))^i \right) J_{M_H (\alpha)\otimes \Omega_{M^r} (\alpha), \underline{h}}(w ;\tau ,\underline{z}).\notag
 \end{align}
 \end{lemma}
 
\paragraph{Proof}
 The proof is by induction on $\ell$, the case $\ell =0$ being clear. Suppose the result holds for all $k$, $0\leq k < \ell$.\ The $n=1$ and $p=1$ case of Lemma \ref{LemmaRecursion} gives
 \begin{align}
 &J_{M_H (\alpha)\otimes \Omega_{M^r} (\alpha), \underline{h}}(a [-1]^\ell w ;\tau ,\underline{z}) \nonumber \\
     &\hspace{5mm}= \tr_{M_H (\alpha)\otimes \Omega_{M^r} (\alpha)} o(a) o(a[-1]^{\ell -1} w) \zeta_1^{h_1 (0)}\cdots \zeta_m^{h_m (0)} q^{L(0)-\mathbf{c}/24} \nonumber \\
 &\hspace{10mm} + (\ell -1)\gby{a,a}E_2 (\tau) J_{M_H (\alpha)\otimes \Omega_{M^r} (\alpha), \underline{h}}(a [-1]^{\ell -2} w ;\tau ,\underline{z}) \nonumber \\
     &\hspace{5mm}=\gby{a,\alpha} J_{M_H (\alpha)\otimes \Omega_{M^r} (\alpha), \underline{h}}(a [-1]^{\ell -1} w ;\tau ,\underline{z})\nonumber \\
 &\hspace{10mm} +(\ell -1) \gby{a,a}E_2 (\tau) J_{M_H (\alpha)\otimes \Omega_{M^r} (\alpha), \underline{h}}(a [-1]^{\ell -2} w ;\tau ,\underline{z}) , \notag
 \end{align}
 where the $E_2 (\tau)$ occur because $h_j (0)a=0$ for all $j$, so that $\tilde{E}_2 (\tau ,0)= E_2 (\tau)$.\ Applying the induction hypothesis on 
 \[
 J_{M_H (\alpha)\otimes \Omega_{M^r} (\alpha), \underline{h}}(a [-1]^{\ell -1} w
 ;\tau ,\underline{z}) \quad \text{and}\quad J_{M_H (\alpha)\otimes \Omega_{M^r} (\alpha), \underline{h}}(a [-1]^{\ell -2} w ;\tau ,\underline{z}),
 \]
 we find
 \begin{align*}
 & J_{M_H (\alpha)\otimes \Omega_{M^r} (\alpha), \underline{h}}(a [-1]^{\ell} w ;\tau ,\underline{z}) \\
     &=\gby{a,\alpha} \left( \sum_{0\leq i \leq (\ell -1) /2} c_{\ell -1,\ell -1-2i} \gby{a,\alpha}^{\ell -1-2i} (\gby{a,a}E_2 (\tau))^i \right) J_{M_H (\alpha)\otimes \Omega_{M^r} (\alpha), \underline{h}}(w ;\tau ,\underline{z}) \\
 &\hspace{10mm} +(\ell -1)\gby{a,a}E_2 (\tau) \left( \sum_{0\leq i \leq (\ell -2) /2} c_{\ell -2,\ell -2-2i} \gby{a,\alpha}^{\ell -2-2i} (\gby{a,a}E_2 (\tau))^i \right)\\
 &\hspace{15mm} \cdot J_{M_H (\alpha)\otimes \Omega_{M^r} (\alpha), \underline{h}}(w ;\tau ,\underline{z})\\
     &=\left( \sum_{0\leq i \leq \ell /2} c_{\ell ,\ell -2i} \gby{a,\alpha}^{\ell -2i} (\gby{a,a}E_2 (\tau))^i \right) J_{M_H (\alpha)\otimes \Omega_{M^r} (\alpha), \underline{h}}(w ;\tau ,\underline{z}),
 \end{align*}
 as desired. The last equality holds since
 \begin{align*}
 &\gby{a,\alpha} \left( \sum_{0\leq i \leq (\ell -1) /2} c_{\ell -1,\ell -1-2i} \gby{a,\alpha}^{\ell -1-2i} (\gby{a,a}E_2 (\tau))^i \right)\\
 &\hspace{5mm} +(\ell -1)\gby{a,a}E_2 (\tau) \left( \sum_{0\leq i \leq (\ell -2) /2} c_{\ell -2,\ell -2-2i} \gby{a,\alpha}^{\ell -2-2i} (\gby{a,a}E_2 (\tau))^i \right)\\
     &\hspace{10mm}=\sum_{0\leq i\leq \ell /2} \left(c_{\ell -1,\ell -1-2i}+(\ell -1)c_{\ell -2,\ell -2i}\right)\gby{a,\alpha}^{\ell -2i} (\gby{a,a}E_2 (\tau))^{i},
 \end{align*}
 so that $c_{\ell ,\ell -2i}:=c_{\ell -1,\ell -1-2i}+(\ell -1)c_{\ell -2,\ell
 -2i}$. $\hfill \Box$\\
 
 Note that
 \begin{align}
 J_{M_H (\alpha)\otimes \Omega_{M^r} (\alpha), \underline{h}}(w ;\tau ,\underline{z})
 &=J_{M_H (\alpha)\otimes \Omega_{M^r} (\alpha), \underline{h}}(1\otimes w ;\tau ,\underline{z})\notag \\
 &=J_{M_H (\alpha) ,\underline{h}}(\textbf{1};\tau ,\underline{z}) \cdot J_{\Omega_{M^r} (\alpha), \underline{h}}(w ;\tau ,\underline{z}), \notag
 \end{align}
 while Equation (\ref{HeisenbergJpartition}) gives
 \begin{equation}
 J_{M_H (\alpha)\otimes \Omega_{M^r} (\alpha), \underline{h}}(w ;\tau ,\underline{z})
=\frac{\zeta_1^{\gby{\alpha ,h_1}}\cdots \zeta_m^{\gby{\alpha ,h_m}} q^{\gby{\alpha ,\alpha}}}{\eta (\tau)^d} J_{\Omega_{M^r} (\alpha), \underline{h}}(w ;\tau ,\underline{z}).
 \label{quasi2}
 \end{equation}
Set
\[
 g_{\ell,i,a}(\tau) := c_{\ell ,\ell -2i} (\gby{a,a}E_2 (\tau))^i 
\]
and
\[
f_{a,\alpha ,\ell}(\tau) :=\! \sum_{0\leq i\leq \ell /2} c_{\ell ,\ell -2i}
\gby{a,\alpha}^{\ell -2i} (\gby{a,a}E_2 (\tau))^i =\!\sum_{0\leq i\leq \ell /2} g_{\ell,i,a}(\tau) \gby{a,\alpha}^{\ell -2i} .
\]
Combining Lemma \ref{quasi1} and (\ref{quasi2}) establishes
\begin{align}
&J_{M_H (\alpha)\otimes \Omega_{M^r} (\alpha), \underline{h}}(a [-1]^\ell w ;\tau ,\underline{z})= f_{a,\alpha ,\ell} (\tau)\frac{\zeta_1^{\gby{\alpha ,h_1}}\cdots \zeta_m^{\gby{\alpha ,h_m}} q^{\gby{\alpha ,\alpha}/2}}{\eta (\tau)^d} J_{\Omega_{M^r} (\alpha), \underline{h}}(w ;\tau ,\underline{z}) .
\notag 
\end{align}
We now take $u_1 ,\dots ,u_d$ to be an orthogonal basis for $H$ and let $\ell_1
,\dots, \ell_d$ be nonnegative integers. We first prove Theorem \ref{TheoremBeta} for elements of the form $v=u_1 [-1]^{\ell_1} \cdots u_d [-1]^{\ell_d} w$. In this case, Lemma \ref{quasi1} implies
\begin{align*}
&J_{M_H (\alpha)\otimes \Omega_{M^r} (\alpha), \underline{h}}(v ;\tau ,\underline{z})\notag \\
&\hspace{5mm}=f_{u_1 ,\alpha ,\ell_1} (\tau)\cdots f_{u_d ,\alpha ,\ell_d} (\tau)\frac{\zeta_1^{\gby{\alpha ,h_1}}\cdots \zeta_m^{\gby{\alpha ,h_m}} q^{\gby{\alpha ,\alpha}/2}}{\eta (\tau)^d} J_{\Omega_{M^r} (\alpha), \underline{h}}(w ;\tau ,\underline{z}) .
\end{align*}
Recalling the module decomposition (\ref{Decomp2}) for $M^r$, it follows that
\begin{equation}
\begin{aligned}
&J_{r, \underline{h}}(v ;\tau ,\underline{z})=\sum_{t=1}^\delta \sum_{\alpha \in \Lambda +\gamma_t} J_{M_H (\alpha)\otimes \Omega_{M^r} (\alpha), \underline{h}}(v ;\tau ,\underline{z}) \\
 &\hspace{0mm}=\sum_{t=1}^\delta \sum_{\alpha \in \Lambda +\gamma_t}
 f_{u_1 ,\alpha ,\ell_1} (\tau)\cdots f_{u_d ,\alpha ,\ell_d} (\tau) \frac{\zeta_1^{\gby{\alpha ,h_1}}\cdots
 \zeta_m^{\gby{\alpha ,h_m}} q^{\gby{\alpha ,\alpha}/2}}{\eta (\tau)^d}
 J_{\Omega_{M^r} (\alpha), \underline{h}}(w ;\tau ,\underline{z}) .
\end{aligned}
\label{quasi4}
\end{equation}
 Each $\alpha \in \Lambda +\gamma_t$ may be written as $\alpha =u+\gamma_t$ for some $u\in \Lambda$. The isomorphism (\ref{OmegaIso}) then shows
 \[
 \Omega_{M^r} (\alpha) = \Omega_{M^r} (u +\gamma_t)\cong \Omega_{M^r} (\gamma_t).
 \]
 Therefore, (\ref{quasi4}) becomes
 \begin{align}
     &J_{r, \underline{h}}(v ;\tau ,\underline{z}) \notag \\
     &=\sum_{t=1}^\delta \frac{J_{\Omega_{M^r} (\gamma_t), \underline{h}}(w ;\tau ,\underline{z})}{\eta (\tau)^d} \sum_{\alpha \in \Lambda +\gamma_t}
 f_{u_1 ,\alpha ,\ell_1} (\tau)\cdots f_{u_d ,\alpha ,\ell_d} (\tau) \zeta_1^{\gby{\alpha ,h_1}}\cdots \zeta_m^{\gby{\alpha ,h_m}} q^{\gby{\alpha ,\alpha}/2} \notag \\
 &=\eta (\tau)^{-d}\sum_{t=1}^\delta J_{\Omega_{M^r} (\gamma_t),
 \underline{h}}(w ;\tau ,\underline{z})\sum_{\alpha \in \Lambda +\gamma_t} \sum_{i_1=0}^{\ell_1 /2} \cdots \sum_{i_d=0}^{\ell_d /2}  g_{\ell_1,i_1,u_1}(\tau) \cdots g_{\ell_d,i_d,u_d}(\tau) \notag \\
  &\quad\ \  \cdot \gby{u_1 ,\alpha}^{\ell_1 -2i_1} \cdots \gby{u_d ,\alpha}^{\ell_d -2i_d}
 \zeta_1^{\gby{\alpha ,h_1}}\cdots \zeta_m^{\gby{\alpha ,h_m}} q^{\gby{\alpha
 ,\alpha}/2}.\label{quasi5}
 \end{align}
 Since the terms $g_{\ell_j ,i_j, u_j}(\tau)$ are independent of $\alpha$, Equation (\ref{quasi5}) becomes

 \begin{align}
     &\eta (\tau)^{-d}\sum_{t=1}^\delta J_{\Omega_{M^r} (\gamma_t), \underline{h}}(w ;\tau ,\underline{z}) \sum_{i_1=0}^{\ell_1 /2} \cdots \sum_{i_d=0}^{\ell_d /2} g_{\ell_1,i_1,u_1}(\tau) \cdots g_{\ell_d,i_d,u_d}(\tau) \notag \\
 &\cdot \sum_{\alpha \in \Lambda +\gamma_t}
 \gby{u_1 ,\alpha}^{\ell_1 -2i_1} \cdots \gby{u_d ,\alpha}^{\ell_d -2i_d}
 \zeta_1^{\gby{\alpha ,h_1}}\cdots \zeta_m^{\gby{\alpha ,h_m}} q^{\gby{\alpha ,\alpha}/2}. \notag
 \end{align}
 Finally, the functions
 \begin{equation}
 \sum_{\alpha \in \Lambda +\gamma_t} \gby{u_1 ,\alpha}^{\ell_1 -2i_1} \cdots \gby{u_d ,\alpha}^{\ell_d -2i_d}
\zeta_1^{\gby{\alpha ,h_1}}\cdots \zeta_m^{\gby{\alpha ,h_m}} q^{\gby{\alpha ,\alpha}/2} \notag
 \end{equation}
 are linear combinations of functions of the form
 \begin{equation}
 \sum_{\alpha \in \Lambda +\gamma_t} \gby{a,\alpha}^{\ell_1 +\cdots +\ell_d -2(i_1 +\cdots +i_d)}\zeta_1^{\gby{\alpha ,h_1}}\cdots \zeta_m^{\gby{\alpha ,h_m}} q^{\gby{\alpha ,\alpha}/2} ,\notag
 \end{equation}
 for various $a \in H$.\ These are the functions (\ref{ThetaIntro}) above and the functions $\theta_{\underline{h}}$ considered in \cite{KrauelI}.\ Finally, we consider arbitrary $v$ as in (\ref{introeq2}).\ Note that if any of the $m_{x,y}$ ($1\leq x\leq d, 1 \leq y \leq \nu_d$) do not equal $1$, then applications of Lemma \ref{LemmaRecursion} will reduce $J_{r, \underline{h}}(v ;\tau ,\underline{z})$ to sums of the form
 \[
 H(\tau) J_{r, \underline{h}}(u_1 [-1]^{\ell_1} \cdots u_d [-1]^{\ell_d} w ;\tau ,\underline{z})
 \]
 for appropriate $\ell_1 ,\dots ,\ell_d \in \mathbb{N}$ and quasi-modular form $H(\tau)$ of weight prescribed in the statement of Theorem \ref{TheoremBeta}. The proof of Theorem \ref{TheoremBeta} is now complete. 
 

%

\begin{thebibliography}{99}
%
%
\bibitem{AM}
 G.~Anderson and G.~Moore,
\newblock \emph{Rationality in conformal field theory}. Comm. Math. Phys., \textbf{117}, 441--450 (1988).

\bibitem{Buhl-Spanning}
 G.~Buhl, \emph{A spanning set for {VOA} modules}. J.~Alg., \textbf{254}(1), 125--151 (2002).

\bibitem{DLM-orbifold}
 C.~Dong, H.~Li and G.~Mason, \emph{Modular-invariance of trace functions in orbifold theory and
              generalized {M}oonshine}. Comm. Math. Phys., \textbf{214}(1), 1--56 (2000).

\bibitem{DLM-twisted}
 C.~Dong, H.~Li and G.~Mason, \emph{Twisted representations of vertex operator algebras and
              associative algebras}. Internat. Math. Res. Notices., 389--397 (1998).

\bibitem{DLMa-elliptic}
 C.~Dong, K.~Liu and X.~Ma,  \emph{Elliptic genus and vertex operator algebras}.
  Pure Appl. Math. Q., \textbf{1}(4, part 3), 791--815 (2005).

\bibitem{DM-Effective}
 C.~Dong and G.~Mason,  \emph{Rational vertex operator algebras and the effective central
              charge}. Int. Math. Res. Not., 2989--3008 (2004).

\bibitem{DM-Shifted}
 C.~Dong and G.~Mason,  \emph{Shifted vertex operator algebras}. Math. Proc. Cambridge Philos. Soc., \textbf{141}(1), 67--80 (2006).

\bibitem{DM-Integrability}
 C.~Dong and G.~Mason,  \emph{Integrability of $C_2$-cofinite Vertex Operator Algebras}. Int. Math. Res. Not. Art.~Id., 80468, 1--15 (2006).

\bibitem{DM-Theta}
 C.~Dong and G.~Mason,  \emph{Transformation laws for theta functions}. Proceedings on {M}oonshine and related topics {M}ontr\'eal, {QC}, Providence, RI, Amer. Math. Soc. (1999);  CRM Proc. Lecture Notes, Vol.~\textbf{30}, 15--26 (2001).
  

\bibitem{DMN-quasi}
C.~Dong,  G.~Mason and K.~Nagatomo, \emph{Quasi-modular forms and trace functions associated to free boson and
  lattice vertex operator algebras}. Internat, Math. Res. Notices, \textbf{8}, 409--427 (2001).
%
\bibitem{EZ}
M.~Eichler and D.~Zagier, \emph{The theory of {J}acobi forms}. Progress in
  Mathematics, Vol.~55, Birkh\"auser Boston Inc., Boston, MA, (1985).

\bibitem{FHL}
I.~Frenkel,  Y.-Z.~Huang and J.~Lepowsky, \emph{On axiomatic approaches to vertex
operator algebras and modules}.  Memoirs Amer. Math. Soc., \textbf{494}, (1993).


\bibitem{GK-differential}
 M.~Gaberdiel and C.~Keller,  \emph{Differential operators for elliptic genera}.  Commun. Number Theory Phys., \textbf{3}(4), 593--618 (2009).

\bibitem{GN}
 M.~Gaberdiel and A.~Neitzke, \emph{Rationality, quasirationality and finite {$W$}-algebras}. Comm. Math. Phys., \textbf{238}(1-2),  305--331 (2003).
 
\bibitem{HE-characters}
 R.~Heluani and J.~Van Ekeren,  \emph{Characters of topological $N=2$ vertex algebras are Jacobi forms on the moduli space of elliptic supercurves}. {\tt arXiv:1405.6128v2} (2014).

\bibitem{Kac-infinite}
 V.~Kac,  \emph{Infinite-dimensional {L}ie algebras}.  2nd Edition, Cambridge University Press,
 Cambridge, (1985).
 
\bibitem{KYY} 
  T.~Kawai,  Y.~Yamada  and S.-K.~Yang, \emph{Elliptic genera and $N=2$ superconformal field theory}. Nucl. Phys. B, \textbf{414}, 191--212 (1994).
 
\bibitem{Kaw}
 T.~Kawai and K.~Yoshioka, \emph{String partition functions and infinite products}. Adv. Theor. Math. Phys., \textbf{4}(2), 397--485 (2000).

\bibitem{KrauelThesis}
M.~Krauel, 
\newblock \emph{Vertex operator algebras and {J}acobi forms} ({P}h{D} {D}issertation). ProQuest Dissertations and Theses (1039264015). ISBN: 978-1-2675-3373-9 (2012).

\bibitem{KrauelMasonI}
M.~Krauel and G.~Mason,
\newblock \emph{Vertex operator algebras and weak {J}acobi forms}. Int. J.~Math., \textbf{23}(6), 1250024--1250034 (2012).

\bibitem{KrauelI}
 M.~Krauel,
\newblock \emph{A {J}acobi theta series and its transformation laws}. Int. J.~Number Theory, \textbf{10}(6), 1343--1354 (2014).

\bibitem{LL}
 J.~Lepowsky and H.~Li, \emph{Introduction to vertex operator algebras and their representations}. In: Progress in Mathematics, \textbf{227},
Birkh\"auser Boston Inc., Boston, MA, (2004).

\bibitem{Li}
 H.~Li, \emph{Local systems of twisted vertex operators, vertex operator superalgebras and twisted modules}. In: Moonshine, the {M}onster, and related topics ({S}outh {H}adley, {MA}, 1994), Contemp. Math., \textbf{193}, Amer. Math. Soc., Providence, RI, 203--236 (1996).

\bibitem{Li-Bilinear}
 H.~Li,  \emph{Symmetric invariant bilinear forms on vertex operator algebras}. J. Pure Appl. Algebra, \textbf{96}(3), 279--297 (1994).

\bibitem{Libgober-Elliptic}
A.~Libgober,  \emph{Elliptic genera, real algebraic varieties and quasi-{J}acobi forms}. In: Topology of stratified spaces,  Math. Sci. Res. Inst. Publ., Vol.~58, Cambridge Univ. Press, Cambridge,  95--120 (2011).
  
\bibitem{Mason-Lattice}
 G.~Mason, 
\newblock \emph{Lattice subalgebras of strongly regular vertex operator algebras}.
 To appear in: Proceedings of the Heidelberg Conference on Vertex Operator Algebras and Related topics, Birkh\"{a}user. {\tt arXiv:1110. 0544v1} (2011).

\bibitem{MTZ}
 G.~Mason, M.~Tuite and A.~Zuevsky,  \emph{Torus {$n$}-point functions for
 {$\mathbb{R}$}-graded vertex operator superalgebras and continuous fermion orbifolds}. Commun. Math. Phys., \textbf{283}(2), 305--342 (2008).
   
\bibitem{Miy}
 M.~Miyamoto, \emph{A modular invariance on the theta functions defined on vertex operator algebras}. Duke Math.~J., \textbf{101}(2),  221--236 (2000).

\bibitem{Oberdiek-Serre}
 G.~Oberdieck, 
\newblock \emph{A {S}erre {D}erivative for even weight {J}acobi Forms}. {\tt arXiv: 1209.5628} (2012).
   
\bibitem{Skor-critical}
 N.~Skoruppa,  \emph{Jacobi forms of critical weight and {W}eil representations}.
 In: Modular forms on {S}chiermonnikoog, Cambridge Univ. Press, Cambridge, 239--266 (2008).

\bibitem{Zhu}
 Y.~Zhu, \emph{Modular invariance of characters of vertex operator algebras}. J.~Amer. Math. Soc., \textbf{9}(1), 237--302 (1996).
\end{thebibliography}
%

%

\end{document}